\input epsf
\def\pf(#1,#2){\langle\langle #1,#2\rangle\rangle}
\def\pf(#1){\langle\langle #1\rangle\rangle}
\def\symb(#1){\langle #1\rangle}
\def\Hom{\mathop{\rm Hom}}
\def\Aut{\mathop{\rm Aut}}
\def\sgn{\mathop{\rm sgn}}

\def\Eq{\mathop{\sl Eq}}

\def\tr{\mathop{\rm tr}}
\def\quadratic{symmetric bilinear }
\def\G{SL(2,K)}
\def\Z{{\bf Z}}
\def\R{{\bf R}}
\def\P{{\bf P}}

\def\TccI{{DJ}}
\def\dotK{{K^*}}
\def\dotKK{K^{*2}}
\def\fund#1{\pi_1(SL(2,#1))}
\def\ov{\overline}
\def\Q{{\bf Q}}
\def\mod{\,{\rm mod}\,}
\def\qed{\hfill{$\diamond$}}
\def\cc{\mathop\tau}
\def\cct{\mathop{\tau_c}}
\def\c1c{\mathop\tau}
\def\w{\mathop{\sl w}}
\def\u{u}
\def\ww{w}
\def\com(#1){#1^{[,]}}
\def\dotK{{K^*}}
\def\dotQ{{\Q^*}}
\def\dotQQ{{\Q}^{*2}}
\def\fund#1{\pi_1(SL(2,#1))}
\def\ov{\overline}
\def\wt{\widetilde}
\def\Q{{\bf Q}}
\def\mod{\,{\rm mod}\,}
\newcount\m
\def\t{\the\m\global\advance\m by 1}
\m=1
\newcount\n
\def\f{\the\n\global\advance\n by 1}
\n=1

\centerline{\bf Tautological characteristic classes II: the Witt class}
\bigskip
\rm
Jan Dymara

Instytut Matematyczny Uniwersytetu Wroc\l awskiego,

pl.~Grunwaldzki 2, 50-384 Wroc\l aw

{\tt jan.dymara@uwr.edu.pl}
\smallskip

Tadeusz Januszkiewicz

Instytut Matematyczny PAN,

ul.~Kopernika 18, 51-617 Wroc\l aw

{\tt tjan@impan.pl}

\footnote{}{Both authors were supported by
Polish NCN grant UMO-2016/23/B/ST1/01556.}

\bigskip\rm

\bigskip
\bigskip
\centerline{\bf Introduction.}
\medskip

Let $H^2(G,U)$ be the second cohomology of a discrete group $G$,
with constant coefficients $U$.
Elements of $H^2(G, U)$
correspond to central extensions of  $G$ by $U$.
If $G$ is perfect, then among all central extensions, with varying $U$,
there is the universal one (cf.~[Brown, Exercise IV.3.7]).
The corresponding class,
which we call the Moore class, lives in
$H^2(G,H_2G)$. 

Cohomology classes of a (discrete) group $G$ give rise to characteristic
classes of 
$G$-bundles and of representations.
Let $\Sigma$ be a closed
oriented surface. 
A principal $G$-bundle 
over $\Sigma$ has a monodromy representation
$\rho\colon \pi_1(\Sigma)\to G$, well-defined up to $G$-conjugation.
Given $\tau\in H^2(G,U)$, one defines the corresponding characteristic class
of the bundle, 
or of $\rho$, by 
$\tau(\rho):=\rho^*\tau\in
H^2(\pi_1(\Sigma),U)\simeq H^2(\Sigma,U)\simeq U$ (the last isomorphism is
given by evaluation on the fundamental class of $\Sigma$).
One important application of these characteristic classes is in the study of
the representation variety $\Hom(\pi_1(\Sigma),G)/G$, i.e.~of the moduli space
of $G$-bundles over $\Sigma$ (cf.~[Lab]). This variety often has additional,
topological or algebro-geometric structure, and one might be interested
in its connected components. The picture one strives to achieve
here is modelled on W.~Goldman's
description of $\Hom(\pi_1(\Sigma),PSL(2,\R))$: for a genus $g$ surface, this space has
$4g-3$ connected components indexed by the Euler class of the representation
(cf.~[Gold88]).
In this classical case there  are also natural continuous deformations
of representations (or flat bundles) called Fenchel--Nielsen twists. 
For an arbitrary group $G$ it is
somewhat vague
what the ``connected components''
of the representation variety are, but the twists do generalize (cf.~[Gold86], or Section 5).
Thus, one wants to consider characteristic classes that are stable under twists.
For a perfect $G$, one may ask whether the universal Moore class is stable
under twists.
It turns out it is not (e.g.~not for $SL(2,K)$). However, a universal twist-stable cohomology
class $\ww_G$ does exist, and there is a beautifully simple condition that we call
\it equicommutativity \rm which is equivalent to twist-stability.
The Moore class coefficient group $H_2G$ has a largest
``equicommutative quotient'' $\Eq(G)$, and
the class $\ww_G\in H^2(G,\Eq(G))$ can be obtained as the image of the
Moore class under the coefficient map $H_2G\to\Eq(G)$.

Let $K$ be an {arbitrary} infinite field.
The above discussion applies to $SL(2,K)$, but here we have more structure.
The action of $SL(2,K)$ 
on its homogeneous space
$\P^1(K)$ gives rise to another class,
defined by J.~Ne\-ko\-v\'a$\check{\rm r}$ in [Ne], 
$\w^I\in H^2(SL(2,K), I^2(K))$, which we call the Witt class. Here $I^2(K)$ is the square of the fundamental ideal
of $W(K)$,
the Witt ring of \quadratic forms over $K$.
It follows from universality that  $\w^I$  is the image of the Moore class
by a certain map $H_2(SL(2,K)) \to I^2(K)$.
We identify this map with
$H_2(SL(2,K))\to \Eq(SL(2,K))$ and prove the following result.

\smallskip
\bf Theorem A \rm (Theorem 10.1)\sl\par
Let $K$ be an infinite field.
The group $\Eq(SL(2,K))$ is isomorphic to $I^2(K)$, and
the Witt class $\w^I\in H^2(SL(2,K),I^2(K))$ is the universal equicommutative
class.
\smallskip\rm

From [\TccI] we know that the Witt class is bounded with respect to the
natural seminorm on $W(K)$.
This is analogous to the classical Milnor--Wood inequality for the Euler class
of flat $SL(2,\R)$-bundles.
Milnor's inequality is sharp: all values allowed by it are indeed Euler
classes of flat bundles.
We study the corresponding saturation problem for the Witt class and prove:

\smallskip
\bf Theorem B \rm(Theorem 11.6)\sl\par
Let $K$ be an infinite field.
\item{(a)} The Witt class of any flat $SL(2,K)$-bundle over an oriented
closed surface of genus $g$ has norm $\le 4(g-1)+2$.
\item{(b)} 
The set of Witt classes of flat $SL(2,K)$-bundles over an oriented
closed surface of genus $g$ contains the set of elements
of $I^2(K)$ of norm $\le 4(g-1)$.
\smallskip\rm
The form of Milnor's inequality we established for
general fields
is not sharp. An example of non-sharpness is constructed over
the field of Laurent series
with rational coefficients.
But for $K=\Q$ we have the sharp result:
\smallskip
\bf Theorem C \rm(Theorem 12.2)\sl\par
The set of Witt classes of all representations of $\pi_1(\Sigma_g)$
in $SL(2,\Q)$
is equal to the set of elements of $I^2\Q$ with norm
$\le4(g-1)$.
\smallskip\rm
To prove Theorem B we use arithmetic properties of Markov
surfaces, established in [GMS],
to construct the required representations.
The proof of Theorem C uses the classical Milnor--Wood inequality and
Meyer's even more classical theorem from the theory of quadratic
forms over $\Q$.
\smallskip
The Witt class can be constructed for $PSL(2,K)$,
but, in general, it is not equicommutative for that group:
that case requires further study.
\smallskip
The paper is divided into four parts, each with its own introduction.
\smallskip
We would like to thank the anonymous referee for numerous
helpful suggestions.

\rm

\bigskip
\centerline{\bf I. Moore and Witt classes.}
\bigskip
In this part we present the main two protagonists: the Witt and Moore classes.
Both are constructed as tautological cohomology classes in the sense of [\TccI].
\medskip
\bf 1. Tautological construction of the Witt class.\rm
\medskip
\def\a{1}
\m=1
\n=1
The Witt class was first defined by
Nekov\'a$\check{\rm r}$ (cf.~[Ne]). It is a cohomology class
$\w\in H^2(SL(2,K),W(K))$, where $K$ is an infinite field and $W(K)$
is the Witt group of \quadratic forms over $K$. A tautological
construction of this class is given in [\TccI, Section 7]. We briefly recall
this construction now. Later, in Section 9, we explain how
Nekov\'a$\check{\rm r}$ modified the class $\w$ to $\w^I\in H^2(SL(2,K), I^2(K))$.
In this section $G=SL(2,K)$.

The group $G$ acts on $\P^1(K)$. The infinite simplex with vertex set
$\P^1(K)$ is a contractible $G$-simplicial complex that we denote by $X$. 
It carries a $G$-invariant tautological cocycle $T$, defined as follows.
First, we define a 2-cochain on $X$ with values in $C_2X$
(the chain group of $X$ with integer coefficients) by assigning
to a 2-simplex this same 2-simplex treated as an element of $C_2X$.
This cochain is not closed---we force it to become closed by
applying to its coefficient group the quotient map $C_2X\to C_2X/B_2X$.
The result is closed (is a cocycle), but it is not $G$-invariant.
We force it to become $G$-invariant by passing to $G$-coinvariants,
i.e.~by applying to it another quotient map
$C_2X/B_2X\to (C_2X/B_2X)_G$. This finally gives $T$. Of course, one fears
that the quotient group $(C_2X/B_2X)_G$ is trivial; however,
quite miraculously, it turns out to be isomorphic to $W(K)$, the additive
group of the Witt ring. (A basic discussion of the Witt ring 
can be found in [EKM, Chapter I].)

The cocycle $T$ can be pulled back to $G$ via (any) orbit map.
In more detail, for any $x\in \P^1(K)$ we consider the $W(K)$-valued
2-cocycle on $G$ defined by
$$(g_0,g_1,g_2)\mapsto T(g_0x,g_1x,g_2x).\leqno(\a.\f)$$
(If $(g_0x,g_1x,g_2x)$ is a degenerate simplex in $X$, the right hand side
is interpreted as zero.) The cohomology class $\w\in H^2(SL(2,K),W(K))$
of this cocycle does not depend on the choice of $x$---this is the Witt class.

To be more explicit we recall the standard Witt group notation:
for $a\in\dotK(:=K\setminus\{0\})$ we denote by $[a]$ the element of $W(K)$ represented
by the 1-dimensional form $ax^2$. The symbol $[0]$ is interpreted as $0$.
The Witt class is represented by the
following (homogeneous) cocycle:
$$(g_0,g_1,g_2)\mapsto [|g_0v,g_1v|\cdot|g_1v,g_2v|\cdot|g_2v,g_0v|];
\leqno(\a.\f)$$
here $v$ is any non-zero vector in $K^2$, and $|g_iv,g_jv|$ stands for
the determinant of the pair of vectors $(g_iv,g_jv)$.
The cocycle depends on $v$, but its cohomology class does not.
In the
non--homogeneous setting we obtain the following cocycle
representing the Witt class:
$$
w(a,b)=
[|v,av|\cdot|av,abv|\cdot|abv,v|]
=[-|v,av|\cdot|v,abv|\cdot|v,bv|].
\leqno(\a.\f)
$$
The standard choice of $v$ is $v=e={1\choose0}$; for this $v$ we get
$$
w(a,b)=[-|e,ae|\cdot|e,abe|\cdot|e,be|]=
[-a_{21}\cdot(ab)_{21}\cdot b_{21}].
\leqno(\a.\f)$$
(In this formula $a_{21}$ denotes the $21$-entry of the matrix $a$.)
This explicit formula will be very useful later.

For details on all the claims made above we again refer the reader
to [\TccI], especially to Section 7 therein.
\medskip
\bf 2. Tautological construction of the Moore class.\rm
\medskip

\def\a{2}
\m=1
\n=1
We will construct a tautological class starting from  the action
of a group $G$ on the standard model of $EG$ (cf.~[Hatcher, Example 1.B.7]).
This model is a $\Delta$-complex with $n$-simplices given by $(n+1)$-tuples 
$[g_0,g_1,\ldots,g_n]$ of elements of $G$. The $G$-action is
$g[g_0,g_1,\ldots,g_n]=[gg_0,gg_1,\ldots,gg_n]$. The quotient, $BG$, has $n$-simplices
given by the orbits of $G$ on the set of $n$-simplices of $EG$. The standard notation is:
$[g_1|g_2|\ldots|g_n]$ for the orbit of $[1,g_1,g_1g_2,\ldots,g_1g_2\ldots g_n]$.
In particular, in $BG$ we have: one vertex $[]$; a loop $[g]$ for each $g\in G$;
a triangle $[g|h]$ for every pair $(g,h)\in G\times G$, with boundary glued to the
edges $[g]$, $[h]$, $[gh]$.

The tautological $n$-cochain for the $G$-action on $EG$ assigns to
an $n$-simplex of $EG$ this same simplex treated as an element of $C_nEG$.
We turn this cochain into a cocycle by dividing the coefficient group by $B_nEG$---this
produces a tautological cocycle in $Z^n(EG,C_nEG/B_nEG)$. We make this cocycle $G$-invariant
by passing to the $G$-coinvariants $U_n$ of the coefficients:
$$U_n:=(C_nEG/B_nEG)_G=(C_nEG)_G/(B_nEG)_G=C_nBG/B_nBG.$$
The resulting $G$-invariant cocycle descends to an element
$T\in Z^n(BG,U_n)$. We also get a cohomology class 
$\tau\in H^n(BG,U_n)\simeq H^2(G,U_n)$.
The inclusion $Z_nBG\to C_nBG$ induces an  inclusion
$Z_nBG/B_nBG\to C_nBG/B_nBG$, thus exhibiting 
$H_nBG=H_nG$ as a subgroup of $U_n$. However,
the values of the cocycle $T$ are 
usually not contained in this subgroup.

Now we specialize to the case of a perfect group $G$ and to $n=2$.
Recall that $G$ is perfect if it has trivial abelianisation, $H_1G=0$.
(We denote by  $H_nG$ the homology group with integer coefficients,
$H_nG=H_n(G,\Z)$.) It follows that $H^1(G,A)=\Hom(G,A)=0$ for all abelian
groups $A$. 

In this special case, the value $T([g|h])$ is the class of $[g|h]$ in $U_2$,
and 
$$\partial[g|h]=[g]-[gh]+[h],$$
which is never zero, so that $T([g|h])\not\in H_2G$ . We will show, however, 
that $T$ is cohomologous to an $H_2G$-valued cocycle.
For each $g\in G$ the edge $[g]$ in $BG$ is a loop; since $G$ is perfect
($H_1G=H_1BG=0$) this loop is null-homologous. Choose, for each $g\in G$, a 2-chain 
$n([g])\in C_2BG$ so that $\partial(n([g]))=[g]$ (we dub $n$ ``the Nekov\'a$\check{\rm r}$
correcting chain''). Then $n\in C^1(G,C_2BG)$, $T$ is cohomologous to $T-\delta n$, and $T-\delta n\in Z^2(G,H_2BG)$:
$$\partial ((T-\delta n)([g|h]))= \partial(T([g|h]))-\partial(n(\partial[g|h]))=
[g]-[gh]+[h]-\partial n([g])+\partial n([gh])-\partial n([h])=0.$$
The coho\-mo\-lo\-gy class $[T-\delta n] \in H^2(G,H_2G)$
is mapped to
the cohomology class $[T] \in H^2(G,C_2BG/B_2BG)$
by the map $\iota\colon H^2(G,H_2G)\to H^2(G,C_2BG/B_2BG)$ (induced by the coefficients inclusion). 
However, there is at most one cohomology class in $H^2(G,H_2G)$ with this property:
indeed, from the Bockstein sequence
$$\ldots\longrightarrow
H^1(G,C_2BG/Z_2BG)\longrightarrow
H^2(G,H_2G)\mathop{\longrightarrow}\limits^{\iota}
H^2(G,C_2BG/B_2BG)\longrightarrow
\ldots
$$
we see that the map $\iota$ is injective (due to $G$ being perfect, $H^1(G,*)=0$).
It follows that the cohomology class of $T-\delta n$ in $H^2(G,H_2G)$ does not depend on the choice of $n$.
\smallskip
\bf Definition \a.\t.\rm\par\nobreak
For a perfect group $G$, the class $\u_G=[T-\delta n]\in H^2(G,H_2G)$ is called
the \it Moore class \rm of $G$.
\smallskip
\bf Remark \a.\t. \rm\par
For perfect $G$ the universal coefficients (evaluation) map
$H^2(G,A)\to\Hom(H_2G,A)$
is an isomorphism; in particular, $H^2(G,H_2G)\simeq\Hom(H_2G,H_2G)$.
The image of the Moore class under this isomorphism is $id_{H_2G}$.
Indeed, for every homology class $[x]\in H_2G$ we have
$$\langle \u_G,[x]\rangle=\langle [T-\delta n],[x]\rangle=\langle T,x\rangle-\langle \delta n,x\rangle=
\langle T,x\rangle=x.$$
This property can serve as another (more standard) definition
of the Moore class---a point of view that will reappear in Section 6.

\bigskip
\centerline{\bf II. Central extensions and characteristic classes.} 
\bigskip
Every cohomology class $\cc{} \in H^2(G,U)$ corresponds to a central extension $\ov{G}$ of $G$ with kernel $U$.
This extension can be used to study $\cc{}$ considered as a characteristic class (i.e.~evaluated 
on $G$-bundles). The first (that we know of) instance of this sort of study 
is Milnor's paper [Mil].
In Section 3 we recall what it means to evaluate $\cc{}$ on a bundle $P$, and how Milnor expressed the result
$\cc(P)$ in term of the lifts to $\ov{G}$ of the monodromies of $P$. In Section 4 we use Milnor's expression to
prove several formulae computing $\cc(P)$ for a bundle $P$ over a surface in terms of restrictions of $P$ to 
subsurfaces. This will be used crucially in Section 11. In Section 5 we discuss twists---natural operations that change bundles.
We compute how $\cc{}$ of a bundle changes under twists, and derive an algebraic condition
on the corresponding central extension $\ov{G}\to G$ that is equivalent to twist--invariance of $\cc{}$.
We call this algebraic condition \it equicommutativity. \rm
In Section 6 we
discuss the Moore class again, this time as the universal class; this allows us
to prove that for perfect $G$ there exists a universal twist--invariant class.

\medskip
\bf 3. Characteristic class in terms of monodromies.\rm
\medskip
\def\a{3}
\m=1
\n=1

In this section $G$ is an arbitrary group, and $U$ is an abelian group.
We consider a class $\cc{} \in H^2(G,U)$. We now recall how this
class can be regarded as a
characteristic class.
Any $G$-bundle $P$ over a space $B$ has a classifying map,
i.e.~a map $B\to BG$ (unique up to homotopy) such that the
pull-back via this map of the universal
$G$-bundle $EG\to BG$ is isomorphic to $P$.
The pull-back of $\cc$ via the classifying map yields
an element $\cc(P)\in H^2(B,U)$---the characteristic class
(corresponding to $\cc$) of the bundle $P$.
We will be interested in the more specific situation when
the base $B$ is 
a closed surface $\Sigma=\Sigma_g$ of
genus $g\ge1$. 
The class $\cc(P)$ can then be evaluated on the fundamental class $[\Sigma]$
to yield an element of $U$.
Evaluation on $[\Sigma]$ defines an isomorphism $H^2(\Sigma,U)\to U$
(e.g.~by the universal coefficient theorem),
so that there is no loss of information in passing from $\cc(P)$ to
$\langle\cc(P),[\Sigma]\rangle$.

The following lemma is well--known (cf.~[Mil]).

\smallskip
\bf Lemma \a.\t.\sl\par
Let $P$ be a $G$-bundle over the surface $\Sigma$, and let $\cc{}\in H^2(G,U)$.
Choose loops $a_1,b_1,\ldots,a_g,b_g$ based at the same point that cut $\Sigma$
into a $4g$-gon and generate $\pi_1(\Sigma)$ with the standard
presentation ($\prod_{i=1}^g[a_i,b_i]=1$). Let $A_i,B_i$ be the monodromy of $P$
along $a_i,b_i$ respectively. Then
$$\langle \cc(P),[\Sigma]\rangle=\prod_{i=1}^g[\ov{A}_i,\ov{B}_i],\leqno(\a.\f)$$
where, for  $g\in G$, we denote by $\ov{g}$ a lift of $g$ to the central extension
$G^{\cc{}}$ of $G$ determined by $\cc$.
\rm\smallskip

The central extension $G^{\cc{}}$ mentioned in the lemma can be described as follows.
Let $\cc$ be represented by a homogeneous cocycle $z\colon G\times G\times G\to U$.
The associated non-homogeneous cocycle is given by $c(g,h)=z(1,g,gh)$. Then, on the set
$G\times U$, we define the multiplication by
$$(g,u)\cdot(g',u')=(gg',uu'\cdot c(g,g')).\leqno(\a.\f)$$
The cocycle condition is equivalent to associativity. We will typically use the standard
(set-theoretic) lift of $G$ to $G^{\cc{}}$: $\ov{g}=(g,1)$.
The abelian group $U$ (=$\{1\}\times U$) is contained in the centre of $G^{\cc{}}$. 
To prove this one checks that $c(g,1)=c(1,g)$ (by setting $g=h$ in the cocycle condition 
$c(g,1)c(g\cdot1,h)=c(g,1\cdot h)c(1,h)$).
We will abbreviate $(1,u)$ to $u$, and use multiplicative notation in $U$. 
(Eventually, for $U=W(K)$, we will switch to the additive convention.)

Conversely, for any central extension $1\to U\to \ov{G}\to G\to 1$
we may choose a set theoretic lift $G\ni g\mapsto\ov{g}\in\ov{G}$ and define
$c\colon G\times G\to U$ by $\ov{g}\cdot\ov{h}=\ov{gh}\cdot c(g,h)$. A change
of the lift changes $c$ within its cohomology class.
We will always assume that the lift of the neutral element of $G$ is
the neutral element of $\ov{G}$ (more obscurely: $\ov{1}=1$).
This assumption implies that $c(1,g)=c(g,1)=1$ for all $g\in G$.

A more thorough discussion of central extensions can be found in [Brown, Chapter IV].

\smallskip
Proof (of Lemma \a.1; cf.~[Mil]). 
Let us first describe the 2-skeleton of the model of
$BG$ that is suitable for us (cf.~[Hatcher, Example 1B.7]).
There is one vertex; for each element $g\in G$ there is an (oriented) edge-loop $[g]$---these form
the 1-skeleton. Then, for every pair $(g,h)\in G^2$, there is a triangle $[g|h]$, with
sides glued to the 1-skeleton along $[g]$, $[h]$ and $[gh]^{-1}$ (going around).
A cocycle $c\in H^2(G,A)\cong H^2(BG,A)$
evaluates on the triangles as $c([g|h])=c(g,h)$.

The surface $\Sigma$ can be expressed as a (convex) polygon $Q$, with $4g$ sides suitably glued in pairs.
We label the vertices of $Q$ by $0,1,\ldots,4g-1$ (counterclockwise) and the edges
(starting by $(0,1)$ and continuing counterclockwise) by $a_1,b_1,a_1^{-1},\ldots,b_{g}^{-1}$.
\midinsert
\epsfxsize=6cm 
\centerline{\magnification1000\epsffile{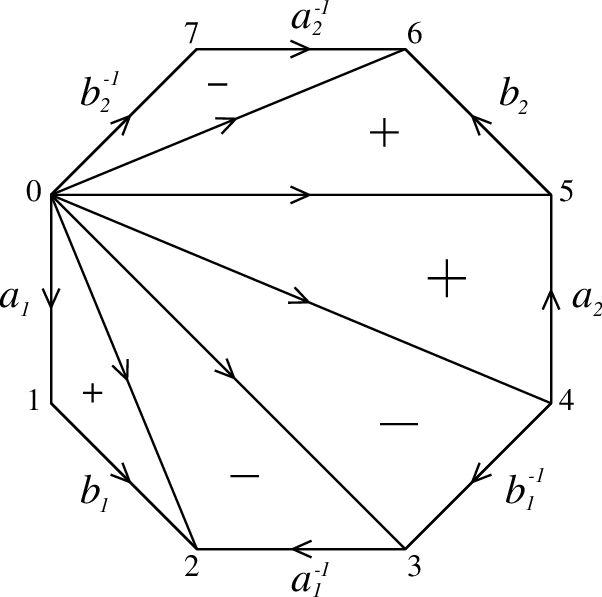}}
\medskip
\centerline{\vbox{\hsize 13cm \noindent
Figure 1: The polygon $Q$ with labels for $g=2$.
The $\Delta$-complex structure
on $\Sigma$ is determined by the arrows.
The fundamental cycle of $\Sigma$ is the
sum of the triangles with the indicated signs.
The map $f$ maps
$\Delta_1=(0,1,2)$ to $[g_1|B_1]=[A_1|B_1]$,
$\Delta_2=(0,3,2)$ to $[g_3|A_1]=[A_1B_1A_1^{-1}|A_1]$, etc.
}}
\endinsert
We also put a $\Delta$-complex structure (cf.~[Hatcher, Section 2.1])
on $\Sigma$.
We divide $Q$ into triangles,
drawing line segments $(0,i)$ for $i=2,3\ldots,4g-2$,
as shown in Figure 1. Then we order the vertices
along each edge, and in each triangle of the triangulation, in a
compatible way. The orders are indicated by arrows in Figure 1.
The arrows are compatible with the boundary gluings, hence we get
a $\Delta$-complex structure on $\Sigma$. (Notice that
the arrows on the boundary edges cannot all be directed counterclockwise
because of the requirement of compatibility with the gluings.)

Now we describe a classifying map $f\colon\Sigma\to BG$ of the bundle $P$.
For convenience, let $c_i^{\epsilon_i}$ be the label of the edge $(i,i+1)$, and let $C_i^{\epsilon_i}$ (equal to
some $A_j^{\pm1}$ or $B_j^{\pm1}$) be the monodromy along that edge.
The map $f$ sends the edge $(i,i+1)$ to $[C_i]^{\epsilon_i}$ (not to $[C_i^{\epsilon_i}]$!).
We map $(0,i)$ to $[g_i]$, where
$g_i=C_0^{\epsilon_0}\cdot\ldots\cdot C_{i-1}^{\epsilon_{i-1}}$. For $\epsilon_i=+1$ we define the triangle
$\Delta_i=(0,i,i+1)$ and map it to $[g_i|C_i]$; for $\epsilon_i=-1$ we define the triangle
$\Delta_i=(0,i+1,i)$ and map it to $[g_{i+1}|C_i]$. Then the fundamental class of $\Sigma$
is represented by the cycle
$\sum_{i=1}^{4g-2}\epsilon_i\Delta_i$,
mapped by $f$ to the cycle
$$\sum_{i\mid\epsilon_i>0}[g_i|C_i]-\sum_{i\mid\epsilon_i<0}[g_{i+1}|C_i]
\leqno(\a.\f)
$$
on which the cocycle $c$ evaluates to
$$\langle \cc(P),[\Sigma]\rangle=\prod_{i\mid\epsilon_i>0}c(g_i,C_i)\prod_{i\mid\epsilon_i<0}c(g_{i+1},C_i)^{-1}.\leqno(\a.\f)$$
On the other hand, $\prod_{i=1}^g[\ov{A}_i,\ov{B}_i]=\prod_{i=0}^{4g-1}\ov{C}_i^{\epsilon_i}$.
Using  the identities
$$\ov{g}\cdot\ov{h}=\ov{gh}\cdot c(g,h),\qquad \ov{g}\cdot\ov{h}^{-1}=\ov{gh^{-1}}\cdot c(gh^{-1},h)^{-1}\leqno(\a.\f)$$
we see that
$$\ov{\prod_{i=0}^{k-1}C_i^{\epsilon_i}}\cdot \ov{C}_k^{\epsilon_k}=\ov{g}_k\cdot\ov{C}_k^{\epsilon_k}=
\cases{
\ov{g}_{k+1}\cdot c(g_k,C_k)=\ov{\prod_{i=0}^kC_i^{\epsilon_i}}\cdot c(g_k,C_k)& if $\epsilon_k=+1$,\cr
\ov{g}_{k+1}\cdot c(g_{k+1},C_k)^{-1}=\ov{\prod_{i=0}^kC_i^{\epsilon_i}}\cdot c(g_{k+1},C_k)^{-1}& if $\epsilon_k=-1.$
}\leqno(\a.\f)$$
We apply this inductively; since $\ov{\prod_{i=0}^{4g-1}C_i^{\epsilon_i}}=\ov{1}=1$, we end up with
$$\prod_{i=0}^{4g-1}\ov{C_i}^{\epsilon_i}=\prod_{i\mid\epsilon_i>0}c(g_i,C_i)\prod_{i\mid\epsilon_i<0}c(g_{i+1},C_i)^{-1}=
\langle \cc(P),[\Sigma]\rangle.\leqno(\a.\f)$$
The lemma is proved.\qed

\medskip

\bf 4. Gluing formulae.\rm
\def\a{4}
\m=1
\n=1
\medskip
Let $1\to U\to\ov{G}\to G\to 1$ be an arbitrary  central group extension.
We choose and fix
a set--theoretic lift $G\to\ov{G}$, to be denoted $g\mapsto\ov{g}$.
We call it the standard lift; we assume that it satisfies $\ov{1}=1$.
(We will sometimes use $\wt{g}$ to denote other, non--standard lifts of $g$.)
We denote by $c$
the corresponding cocycle (so that $\ov{g}\cdot\ov{h}=\ov{gh}\cdot c(g,h)$)
and by $\cct{}$ its cohomology class ($\cct{}\in H^2(G,U)$). 

Suppose that $\xi$ is a flat $G$-bundle over an oriented compact surface
$S$ of genus $g$
with one boundary component. A \it boundary framing \rm of $\xi$ will mean 
the following collection of data: a point $s\in\partial S$; a trivialization of the fibre
$\xi_s$; an element $\wt{W}\in\ov{G}$ that lifts the monodromy $W\in G$  of $\xi$ along
$\partial S$ (oriented compatibly with the orientation of $S$ and based at $s$).
We will often abusively say ``boundary framing $\wt{W}$'', because $\wt{W}$ is the part
of the framing data that presumes the rest of it and appears explicitly in many formulas.
Usually we will use $\wt{W}=\ov{W}$---the standard lift of $W$---and then call our framing \it a standard framing\rm.

A \it standard loop collection, \rm i.e.~a collection of loops
$(x_i,y_i\mid i=1,\ldots,g)$, based at $s$, cutting
$\Sigma$ into a $(4g+1)$-gon, generating
$\pi_1(S)$ and satisfying $[x_1,y_1]\cdot[x_2,y_2]\cdot\ldots\cdot[x_g,y_g]=\partial S$
can be chosen (in many ways). Let $X_i,Y_i\in G$ be the monodromies of $\xi$ along these loops.
\smallskip

\bf Definition \a.\t.\rm\par
The \it relative class \rm  $\ov{c}(\xi,\wt{W})$ of the bundle $\xi$ with boundary framing $\wt{W}$
is defined as the element of $U$ given by
$$\ov{c}(\xi,\wt{W}):=[\ov{X}_1,\ov{Y}_1]\cdot[\ov{X}_2,\ov{Y}_2]\cdot
\ldots\cdot[\ov{X}_g,\ov{Y}_g]\cdot\wt{W}^{-1}.\leqno(\a.\f)$$
If the framing is standard (i.e.~$\wt{W}=\ov{W}$), then we put $\ov{c}(\xi)=\ov{c}(\xi,\ov{W})$.
(The bar over $c$ is a reminder that the class is relative, and some framing is presumed.)
\smallskip
\bf Remark \a.\t.\rm\par \nobreak
The name ``relative class'' is intentionally provocative. We expect that there exists a relative characteristic
class $\ov{\tau}_c$ related to $\cct{}$ such that $\ov{\tau}_c(\xi)$ evaluated
on the relative fundamental class of the base of $\xi$ equals $\ov{c}(\xi)$.
\smallskip
\bf Remark \a.\t.\rm\par
Notice that since $\prod[x_i,y_i]=\partial S$, we have $\prod[X_i,Y_i]=W$;
it follows that $\ov{c}(\xi)\in U$. The commutator
$[\ov{X}_i,\ov{Y}_i]$
does not depend on the choice of the lifts
of $X_i$, $Y_i$ to $\ov{G}$: $[\wt{X}_i,\wt{Y}_i]=[\ov{X}_i,\ov{Y}_i]$. We will use different lifts
to our advantage. On the other hand, changing $\wt{W}$ to a different lift $\wt{W}u$ of $W$
results in a change:
$$\ov{c}(\xi,\wt{W}u)=\ov{c}(\xi,\wt{W})\cdot u^{-1}.\leqno(\a.\f)$$
It is not a priori clear that
$\ov{c}(\xi)$ does not depend on the choice of the collection of loops $(x_i,y_i)$;
we will check this shortly.
\smallskip
If a closed surface is cut along a separating simple loop into two pieces,
a $G$-bundle over that surface decomposes into two bundles over the pieces.
The next lemma describes the relation between the (relative) classes of the
three bundles.
\smallskip
\bf Lemma \a.\t.\sl\par\nobreak
Let $\xi$, $\xi'$ be flat $G$-bundles over oriented surfaces $S$, $S'$ with isomorphic
boundary framings $\wt{W}$. The isomorphism of boundary framings allows one to glue
the bundles $\xi$, $\xi'$; the result is a bundle
$\xi\cup\xi'$ over $\Sigma=S\cup_\partial S'$. We orient $\Sigma$ compatibly with $S$
and opposite to $S'$. Then
$$\langle\cct(\xi\cup \xi'),[\Sigma]\rangle=\ov{c}(\xi,\wt{W})\ov{c}(\xi',\wt{W})^{-1}.\leqno(\a.\f)$$
\smallskip\rm
Proof. Let $(x_i,y_i)$ be a standard collection of loops for $S$, $(x_i',y_i')$ one for $S'$.
Then these collections together, in the order $(x_1,y_1,\ldots,x_g,y_g,y'_{g'},x'_{g'},\ldots,y'_1,x'_1)$,
give a standard set of generators of $\pi_1(S\cup_\partial S')$. We have
$$
\eqalign{
\langle\cct(\xi\cup\xi'),[\Sigma]\rangle&=
[\ov{X}_1,\ov{Y}_1]\ldots[\ov{X}_g,\ov{Y}_g]
[\ov{Y'}_{g'},\ov{X'}_{g'}]\ldots[\ov{Y'}_1,\ov{X'}_1]\cr
&=
[\ov{X}_1,\ov{Y}_1]\ldots[\ov{X}_g,\ov{Y}_g]\cdot\wt{W}^{-1}\cdot
\wt{W}\cdot[\ov{X'}_{g'},\ov{Y'}_{g'}]^{-1}\ldots[\ov{X'}_1,\ov{Y'}_1]^{-1}\cr
&=
[\ov{X}_1,\ov{Y}_1]\ldots[\ov{X}_g,\ov{Y}_g]\cdot\wt{W}^{-1}\cdot
\left([\ov{X'}_1,\ov{Y'}_1]\ldots[\ov{X'}_{g'},\ov{Y'}_{g'}]\cdot\wt{W}^{-1}\right)^{-1}\cr
&=
\ov{c}(\xi,\wt{W})\ov{c}(\xi',\wt{W})^{-1}.}
\leqno(\a.\f)$$

\qed
\smallskip
\bf Corollary \a.\t.\sl\par
The relative class $\ov{c}(\xi,\wt{W})$ does not depend on the choice of a standard loop collection $(x_i,y_i)$.
\smallskip\rm
Proof. Choose any $(\xi',S')$ with the same boundary framing as $(\xi,S)$
(it can be just another copy of $(\xi,S)$). Then compute $\cct(\xi\cup\xi')$
as in Lemma \a.4. We get $\cct(\xi\cup\xi')=\ov{c}(\xi,\wt{W})\ov{c}(\xi',\wt{W})^{-1}$
regardless of the choices of the collections $(x_i,y_i)$ and $(x'_i,y'_i)$.
Varying one of these collections while keeping the other fixed we see the claimed independence.
\qed
\smallskip
Now we set the notation for Lemma \a.6.
Let $\xi$ be a bundle over $S$ with a standard boundary framing,
and let $A\in G$.
Then we can change (twist) the framing by $A$. This means that we change the
trivialization $\xi_s\to G$ by $A$; then all monodromies $M$ change to $^AM:=AMA^{-1}$.
In particular, $W$ changes to $^AW$, and the standard framing of the $A$-twisted bundle
is $\ov{^AW}$. We denote by $^A\xi$ the twisted bundle with this framing.
Another natural choice of framing is $^A\ov{W}$. (Formally, we should use
$^{\ov{A}}\ov{W}=\ov{A}\cdot\ov{W}\cdot\ov{A}^{-1}$, but this does not depend on the
choice of the lift $\ov{A}$ of $A$ and will usually be abbreviated to $^A\ov{W}$.)
In general $^A\ov{W}\ne\ov{^AW}$, so that we have two  natural $A$-twisted bundles with
boundary framing:
$(\xi,^A\ov{W})$ and
$^A\xi=(\xi,\ov{^AW})$. 
\smallskip
\bf Lemma \a.\t.\sl\par
We have $$\ov{c}(\xi,^A\ov{W})=\ov{c}(\xi,\ov{W})=\ov{c}(\xi),\qquad
\ov{c}(^A\xi)=\ov{c}(\xi,\ov{^AW})=\ov{c}(\xi){c(A,W)c(^AW,A)^{-1}}.\leqno(\a.\f)$$
\smallskip\rm
Proof.
Let $(x_i,y_i)$ be a standard collection of loops in $S$, and let $(X_i,Y_i)$ be the monodromies of $\xi$
along these loops. The monodromies for the $A$-twisted $\xi$ are $(^AX_i,^AY_i)$. Since the commutator
of lifts does not depend on the choice of the lifts, we have
$[\ov{^AX_i},\ov{^AY_i}]=[^A\ov{X_i},^A\ov{Y_i}]$.
Therefore
$$\eqalign{
\ov{c}(\xi,{}^A\ov{W})
&=\left(\prod[\ov{^AX_i},\ov{^AY_i}]\right)\left(^A\ov{W}\right)^{-1}
=\left(\prod{}^A[\ov{X_i},\ov{Y_i}]\right) {}^A\ov{W}^{-1}\cr
&={}^A\left(\prod[\ov{X_i},\ov{Y_i}]\cdot\ov{W}^{-1}\right)
={}^A\ov{c}(\xi)=\ov{c}(\xi).}\leqno(\a.\f)$$
(Conjugation does not change the central element $\ov{c}(\xi)$.)
We compare this with $\ov{c}(^A\xi)$ using the standard lift $\ov{A}$ of $A$:
$$^A\ov{W}=\ov{A}\cdot\ov{W}\cdot\ov{A}^{-1}=\ov{AW}c(A,W)\ov{A}^{-1}=
\ov{AWA^{-1}}c(AWA^{-1},A)^{-1}c(A,W)=\ov{^AW}{c(A,W)c(^AW,A)^{-1}}.\leqno(\a.\f)
$$
Now the second formula follows from (\a.2) and (\a.7):
$$
\ov{c}(\xi,\ov{^AW})=\ov{c}(\xi,{}^A\ov{W} c(A,W)^{-1} c(^AW,A))=
\ov{c}(\xi,{}^A\ov{W}) c(A,W) c(^AW,A)^{-1}=
\ov{c}(\xi){c(A,W) c(^AW,A)^{-1}}
\leqno(\a.\f)$$
\qed
\smallskip
Now we set up notation for Lemma \a.7 and give another gluing construction, called ``boundary connected sum''.
Let $\xi$, $\xi'$ be bundles with standard boundary framings $\ov{W}$, $\ov{W'}$ over $S$, $S'$.
We glue $s$ to $s'$ as well as the fibres $\xi_s$, $\xi'_{s'}$
(via the framing trivializations). Then we glue a triangle
$\Delta$ to $\partial S\vee\partial S'$, one side along $\partial S$, one along $\partial S'$
(all vertices to $s=s'$), so that the path $(\partial S)(\partial S')$ is homotopic (through $\Delta$)
to the third side. This third side forms the boundary of the obtained surface $\Sigma$ of genus $g+g'$.
The bundle naturally extends to a bundle $\xi\vee\xi'$ over $\Sigma$ (all monodromies are already visible
in $\xi$ and $\xi'$). The trivializations at $s=s'$ agree and trivialize the new fibre at this point.
The boundary monodromy of $\xi\vee\xi'$ is $WW'$; we use the standard boundary framing, with lift $\ov{WW'}$.

\smallskip
\bf Lemma \a.\t.\sl\par
In the above situation,
$$\ov{c}(\xi\vee\xi')=\ov{c}(\xi\vee\xi',\ov{WW'} )=\ov{c}(\xi)\ov{c}(\xi')c(W,W').\leqno(\a.\f)$$
\rm
Proof. Let $(x_i,y_i)$ and $(x'_j,y'_j)$ be standard loop collections in $S$, $S'$.
Together they form a standard loop collection in $\Sigma$, so that
$$
\eqalign{\ov{c}(\xi\vee\xi',\ov{WW'})&=
\prod[\ov{X}_i,\ov{Y}_i]\prod[\ov{X'}_j,\ov{Y'}_j]\cdot \ov{WW'}^{-1}
=\ov{W}\ov{c}(\xi)\ov{W'}\,\ov{c}(\xi')\ov{WW'}^{-1}\cr
&=\ov{W}\cdot
\ov{W'}\cdot
\ov{WW'}^{-1}\cdot
\ov{c}(\xi)
\ov{c}(\xi')
=
c(W,W')\ov{c}(\xi)
\ov{c}(\xi')
}\leqno(\a.\f)
$$
\qed

One last piece of general calculation:
\smallskip
\bf Lemma \a.\t.\sl\par
Let $\xi$ be a bundle over a surface with one boundary component
and genus $1$ with standard
boundary framing, standard loop generators $(x,y)$  and monodromies $(X,Y)$
(with $[X,Y]=W$). Then
$$\ov{c}(\xi)={c(X,Y) c(Y,X)^{-1}c(W,YX)^{-1}}.\leqno(\a.\f)$$
\rm
Proof.
$$\eqalign{
\ov{c}(\xi)&=
[\ov{X},\ov{Y}]\ov{W}^{-1}=\ov{X}\cdot\ov{Y}\cdot(\ov{Y}\cdot\ov{X})^{-1}\ov{W}^{-1}
=\ov{XY}c(X,Y)(\ov{YX}c(Y,X))^{-1}\ov{W}^{-1}\cr
&=\ov{XY}\cdot\ov{YX}^{-1}\ov{W}^{-1}{c(X,Y) c(Y,X)^{-1}}\cr
&=\ov{XY(YX)^{-1}}c(XY(YX)^{-1},YX)^{-1}\ov{W}^{-1}{c(X,Y) c(Y,X)^{-1}}\cr
&=\ov{W}\cdot\ov{W}^{-1}{c(X,Y) c(Y,X)^{-1}c(W,YX)^{-1}}
={c(X,Y) c(Y,X)^{-1}c(W,YX)^{-1}}.
}\leqno(\a.\f)$$
\qed

\medskip
\bf 5. Twists and equicommutativity.\rm
\def\a{5}
\m=1
\n=1
\medskip

In this section we discuss the twist deformations of flat bundles
over surfaces.
These twists are associated to the names of Fenchel and Nielsen
in the Teichm\"uller case (cf.~[Wol]), and to Goldman in
the Lie group case (cf.~[Gold86]).

Let $P$ be a (flat) $G$-bundle over an oriented
surface $\Sigma$. Choose an oriented simple loop $\ell$ in $\Sigma$,
and a base-point $b\in \ell$.
Trivialize $P_b$ (by a right-$G$-equivariant isomorphism $P_b\to G$).
Then, there is a well-defined element $L\in G$ representing
the monodromy of $P$ along $\ell$.
Choose any $V\in Z_G(L)$ (the centralizer of $L$ in $G$),
and trivialize the bundle $P$ along $\ell$ (with ambiguity $L$ at $b$). 
Then cut $\Sigma$ and $P$ along $\ell$, and glue it back by (left)
multiplication by $V$. (To be precise, we define the right-hand side
and the left-hand side of a tubular neighbourhood of $\ell$ in $\Sigma$
using orientations. Then, after cutting the bundle, each
element $p$ of a trivialized fibre $P_x$ at a point $x\in\ell$ is split into
a left-right pair $p_L$, $p_R$. We glue $p_R$ to $Vp_L$.
Since the trivialization along $\ell$ is $L$-ambiguous at $b$, the gluing over $b$
is well-defined only for $V\in Z_G(L)$.)
The result is a new (flat) $G$-bundle $P_{\ell,V}$ over $\Sigma$---
the \it twist \rm of
$P$ by $V$ along $\ell$. 

It is possible to phrase the above definition in a slightly more
invariant way. Suppose we refrain from choosing a trivialization of $P_b$.
Then we still have the monodromy along $\ell$. It is an element
$L$ in $\Aut(P_b)$, the automorphism group of the right $G$-space $P_b$.
(This $\Aut(P_b)$ is non-canonically isomorphic to $G$; possible
isomorphisms arise from trivializations of $P_b$.)
Then for any $V\in Z_{\Aut(P_b)}(L)$ the bundle $P_{\ell,V}$ is well-defined.

Now suppose that we have a central group extension
$1\to U\to\ov{G}\to G\to 1$,
as in the previous section
(with a lift $g\mapsto\ov{g}$, cocycle $c$,
cohomology class $\cc{}\in H^2(G,U)$).

\smallskip
\bf Theorem \a.\t.\sl\par
Let $P$ be a $G$-bundle over a closed oriented surface $\Sigma$.
Let $\ell$ be an oriented simple loop in $\Sigma$, based at $b$.
Let $L\in G$ be the monodromy of $P$ along $\ell$
(with respect to some trivialization of $P_b$), let
$V\in Z_G(L)$,
and let $P_{\ell,V}$ be the twist of $P$. Let $\cc\in H^2(G,U)$
be a cohomology class represented by a cocycle $c$.
Then
$$\langle\cc(P_{\ell,V}),\Sigma\rangle=\langle\cc(P),[\Sigma]\rangle{c(V,L)  c(L,V)^{-1}}.\leqno(\a.\f)$$
\rm\par
Proof. 
We use the commutator product expression from Lemma 3.1. 
The basic calculation (valid for any two commuting elements $V,L\in G$) is
$$
\ov{V}\,\ov{L}=\ov{VL}c(V,L)=\ov{LV}c(V,L)=
\ov{L}\,\ov{V}c(L,V)^{-1}c(V,L)=
\ov{L}\,\ov{V}{c(V,L) c(L,V)^{-1}}
.\leqno(\a.\f)$$
\par
Case 1. The loop $\ell$ does not separate $\Sigma$.
Then the standard presentation
loops $a_1,\ldots,b_g$ in $\Sigma$ can be chosen so that $b_1=\ell$.
If $A_1,\ldots,B_g$
are the elements of $G$ representing the monodromies of the bundle
$P$ along $a_1,\ldots,b_g$ (with $B_1=L$),
then the monodromies of $P_{\ell,V}$ along these loops are represented
by the same elements except
for one change: $A_1$ gets replaced by $A_1V$.
In the commutator product expression the first term $[\ov{A}_1,\ov{B}_1]$ ($=[\ov{A}_1,\ov{L}]$) changes to
$[\ov{A_1V},\ov{B}_1]$ ($=[\ov{A_1V},\ov{L}]$).
Using (\a.2) we compute:
$$\eqalign{[\ov{A_1V},\ov{L}]
&=[\ov{A}_1\ov{V}c(A_1,V)^{-1},\ov{L}]
=[\ov{A}_1\ov{V},\ov{L}]\cr
&=\ov{A}_1\ov{V}\,\ov{L}\,\ov{V}^{-1}\ov{A}_1^{-1}\ov{L}^{-1}=
\ov{A}_1\ov{L}\,\ov{V}{c(V,L) c(L,V)^{-1}}\ov{V}^{-1}\ov{A}_1^{-1}\ov{L}^{-1}
\cr
&=[\ov{A}_1,\ov{L}]{c(V,L) c(L,V)^{-1}}.}\leqno(\a.\f)$$
The claim follows.
\par
Case 2. 
The loop $\ell$ separates $\Sigma$. Then we cut $\Sigma$ and $P$ along
$\ell$ into two components, say $P_0$ over $\Sigma_0$ and $P_1$ over $\Sigma_1$.
The assumptions of the theorem induce (isomorphic) boundary framings for $P_0$ and $P_1$ (except for
lifts of the boundary monodromy $L$---we take standard lifts).
We have $P=P_0\cup P_1$, $P_{\ell,V}={}^V\!P_0\cup P_1$. Lemmas 4.4, 4.6 give
$$\langle\tau(P_{\ell,V}),[\Sigma]\rangle=\ov{c}({}^V\!P_0)\ov{c}(P_1)^{-1}=
\ov{c}(P_0){c(V,L) c({}^V\!L,V)^{-1}}\ov{c}(P_1)^{-1}
=\langle\tau(P),[\Sigma]\rangle{c(V,L) c(L,V)^{-1}}.\leqno(\a.\f)$$
The last equality uses the fact that ${}^V\!L=L$, a consequence of $V\in Z_G(L)$.
\qed

\smallskip
This theorem leads to the following definition.
\smallskip

\smallskip
\bf Definition \a.\t. \rm
A cocycle $c\colon G\times G\to U$ is called \it equicommutative, \rm
if it satisfies $c(g,h)=c(h,g)$ whenever $gh=hg$.
A cohomology class is \it equicommutative \rm if it is represented by an
equicommutative cocycle.
\smallskip
\bf Proposition \a.\t.\sl\par
Let $c$ be a cocycle with cohomology class $\cct{}\in H^2(G,U)$, let
$U\to\ov{G}\to G$ be the corresponding central extension, and let
$g\mapsto\ov{g}$ be the lift corresponding to $c$.
The following conditions are equivalent:
\item{a)} for every commuting pair $g,h\in G$,
the lifts $\ov{g},\ov{h}\in\ov{G}$ commute;
\item{b)} the cocycle $c$ is equicommutative;
\item{c)} the cohomology class $\cct$ is equicommutative;
\item{d)} for every commuting pair $g,h\in G$, every lift of $g$ commutes
with every lift of $h$
\item{e)} every cocycle representing $\cct$ is equicommutative;
\item{f)} for every commutative subgroup $H<G$ the pre-image of $H$
in $\ov{G}$ is commutative.
\rm\smallskip
Proof. It is straightforward to see that:
\item{-} the weak conditions a), b) are equivalent;
\item{-} the strong conditions d), e), f) are equivalent;
\item{-} the strong conditions imply the weak conditions;
\item{-} b) implies c).

To finish the proof we show that c) implies e).
It is enough to check that all 2-coboundaries are equicommutative.
Let $g,h\in G$ be commuting elements, and let $n\in C^1(G,U)$. Then
$$(\delta n)(g,h)=n(h)n(gh)^{-1}n(g)=n(g)n(hg)^{-1}n(h)=(\delta n)(h,g).\leqno(\a.\f)$$
\qed

\smallskip
\bf Definition \a.\t.\rm\par
A cohomology class $\cc{}\in H^2(G,U)$ is \it twist--invariant, \rm
if for every $G$-bundle $P$ over a closed oriented surface $\Sigma$,
and for every twist $P_{\ell,V}$ of that bundle, we have
$\cc(P_{\ell,V})=\cc(P)$.
\smallskip
In Definition \a.4 one could, equivalently, use the condition
$\langle\cc(P_{\ell,V}),[\Sigma]\rangle=\langle\cc(P),[\Sigma]\rangle$.
This is because evaluation on $[\Sigma]$ is an isomorphism
$H^2(\Sigma,U)\to U$.
\smallskip
\bf Corollary \a.\t.\sl\par
A cohomology class is twist--invariant if and only if it is equicommutative.
\rm\par
Proof.
It follows from Theorem \a.1 that equicommutativity implies twist--invariance.
For the converse,
let $c$ be a cocycle representing a twist--invariant class
$\cc{}\in H^2(G,U)$.
For any pair of commuting elements $g,h\in G$
there exists a $G$-bundle $\xi_{(g,h)}$ over $\Sigma_1=T^2$ with
monodromies (along standard generating loops $a_1,b_2$) 
equal to $g$, $h$. This bundle is a twist of $\xi_{(1,h)}$, and 
Theorem \a.1 gives:
$$\cc(\xi_{(g,h)})=\cc(\xi_{(1,h)}){c(g,h) c(h,g)^{-1}}.\leqno(\a.\f)$$
Now the assumption of twist--invariance implies that $c(g,h)=c(h,g)$.
\qed
\smallskip\rm
\bf Corollary \a.\t.\sl\par
The Witt class is twist--invariant.
\rm\par
Proof. We check that the cocycle (1.4) is equicommutative.
Let $a,b\in SL(2,K)$, and suppose that $ab=ba$. Then
$$w(a,b)=[-a_{21}(ab)_{21}b_{21}]=[-a_{21}(ba)_{21}b_{21}]=w(b,a).\leqno(\a.\f)$$
\qed
\smallskip
\bf Remark \a.\t.\rm\par
Unlike the Witt class, the Moore class, in general, is not twist--invariant.
This is more fully explained in Section 8.

\medskip
\bf 6. Universal classes.\rm
\medskip
\def\a{6}
\m=1
\n=1

In this section we specialize our discussion of characteristic classes to perfect groups.
Recall that a group $G$ is perfect if it is equal to its commutator subgroup $[G,G]$.
In homological terms this means that $H_1G=0$ (we denote $H_i(G,\Z)$ by $H_iG$). Consequently,
the universal coefficients map $H^2(G,A)\to\Hom(H_2G,A)$ is an isomorphism for every abelian
group $A$. In particular, for $A=H_2G$ we have a well-defined class $\u_G\in H^2(G,H_2G)$
that corresponds to $id_{H_2G}$ under this isomorphism. 
(Thus, by Remark 2.2, $\u_G$ coincides with the class constructed in Section 2.)
The class $\u_G$ is universal in the following sense:
for any abelian group $A$ and any class $v\in H^2(G,A)$ there exists a unique
homomorphism $f\colon H_2G\to A$ such that $f_*\u_G=v$
(cf.~[Brown, Exercise IV.3.7]).
Also the central extension $1\to H_2G\to\ov{G}\to G\to 1$ defined by $\u_G$ is universal.
(This extension was one of the early reasons for considering the
second homology of a group;
whence the name ``Schur multiplier'' for $H_2G$.) 
It is known that all universal central extensions of a perfect group $G$ are canonically
isomorphic. Quite often the construction of such an extension and the study
of its kernel is the way to calculate $H_2G$ and to describe $\u_G$.
We call the class $\u_G$ the Moore class, because it was
investigated by Moore for $G=SL(2,K)$.
\smallskip
\bf Definition \a.\t.\rm\par\nobreak
Let $G$ be a perfect group, and let $u$ be a cocycle representing the universal class
$\u_G\in H^2(G,H_2G)$. Let $\com(G)$ be the subgroup of $H_2G$
generated by the set
$$\{{u(g,h) u(h,g)^{-1}}\mid g,h\in G, gh=hg\}.\leqno(\a.\f)$$
We put $\Eq(G):=H_2G/\com(G)$. Let $\ww_G\in H^2(G,\Eq(G))$ be the image of the universal class $\u_G$
by the quotient map $q\colon H_2G\to \Eq(G)$.
\smallskip
\bf Remark \a.\t.\rm\par
The set (\a.1) does not depend on the choice of the cocycle $u$.
Indeed, it can be described as the set of commutators, in $\ov{G}$,
of lifts of pairs of commuting elements $g,h\in G$,
as the following calculation shows:
$${u(g,h) u(h,g)^{-1}}={\ov{gh}^{-1}\ov{g}\,\ov{h}}\left(\ov{hg}^{-1}\ov{h}\,\ov{g}\right)^{-1}=
{\ov{gh}^{-1}\ov{g}\,\ov{h}\,\ov{g}^{-1}\ov{h}^{-1}\ov{hg}}=
{\ov{gh}^{-1}[\ov{g},\ov{h}]\ov{hg}}
=[\ov{g},\ov{h}]
.\leqno(\a.\f)$$
The last equality follows from the fact that for commuting $g,h$ the
commutator $[\ov{g},\ov{h}]$ is central in $\ov{G}$. This commutator
does not depend on the choice of lifts of $g,h$, because all possible
lifts are of the form $\ov{g}u$, $\ov{h}v$ with $u,v$ central in $\ov{G}$.
\smallskip
\bf Remark \a.\t.\rm\par
Another, more topological description of the set (6.1): it consists
of ``genus 1 classes'', i.e.~the classes in $H_2G$ that are images of
the fundamental class of the 2-dimensional torus $T^2$ under some map
$T^2\to BG$. Indeed, such a map associates to the generators
of $\pi_1(T^2)$ an (arbitrary) commuting pair $g,h\in G$; the image
of the fundamental class of $T^2$ is then $u(g,h)u(h,g)^{-1}$ by the computation in the proof of Lemma 3.1.
\smallskip

\bf Theorem \a.\t.\sl\par
Let $G$ be a perfect group.
The class $\ww_G\in H^2(G,\Eq(G))$ is a universal equicommutative class
in the following sense:
for every equicommutative cohomology class $v\in H^2(G,A)$ there exists a unique
homomorphism $g\colon \Eq(G)\to A$ such that $g_*\ww_G=v$.

\rm\par
Proof. Let $f\colon H_2G\to A$ be the homomorphism that maps $\u_G$ to $v$.
Choose a cocycle $u$ representing $\u_G$.
Then $f_*u$ is a cocycle representing $v$, hence it is equicommutative.
It follows that, for every commuting pair $g,h\in G$, we have
$$f\left({u(g,h) u(h,g)^{-1}}\right)={f_*u(g,h)\cdot f_*u(h,g)^{-1}}=1;\leqno(\a.\f)$$
therefore $f$ factors through the quotient map $q\colon H_2G\to \Eq(G)$,
i.e.~$f=g\circ q$ for some $g\colon \Eq(G)\to A$. We get
$$g_*\ww_G=g_*q_*\u_G=f_*\u_G=v.\leqno(\a.\f)$$
The uniqueness statement is proved by contradiction. Suppose two different
homomorphisms $g,g'\colon \Eq(G)\to A$ map $\ww_G$ to $v$; then
$g\circ q,g'\circ q\colon H_2G\to A$ are different,
and both map $\u_G$ to $v$---contradiction.
\qed
\smallskip
\bf Remark \a.\t.\rm\par \nobreak
The class $\ww_G\in H^2(G,\Eq(G))$ is, up to a unique isomorphism, the unique
universal equicommutative class; this is a standard consequence
of universality.
\medskip

We finish this section by indicating a more general point of view on universality.
It will not be used later in this paper.
\smallskip
\bf Proposition \a.\t. \sl\par
Let $G$ be a perfect group, $\u_G\in H^2(G,H_2G)$ its Moore class, and let
$\varphi\colon H_2G\to Q$ be a group epimorphism (coefficient reduction map). 
We set $\u_{G,Q}=\varphi_*\u_G$. Suppose that a cohomology class $v\in H^2(G,A)$ satisfies
the following condition: $v(x)=0$ for all $x\in\ker{\varphi}$.
Then there exists a unique group homomorphism $\psi\colon Q\to A$ such that $\psi_*\u_{G,Q}=v$.

\rm\par
Proof. Let $\Psi\colon H_2G\to A$ be the unique map giving $\Psi_*\u_G=v$.
For each $x\in \ker{\varphi}$ we have:
$$\Psi(x)=\Psi(\u_G(x))=\Psi_*\u_G(x)=v(x)=0.$$
It follows that there exists a $\psi\colon Q\to A$
such that $\psi\circ\varphi=\Psi$. Then
$$\psi_*\u_{G,Q}=\psi_*\varphi_*\u_G=\Psi_*\u_G=v.$$
Now we show that the homomorphism $\psi$ is unique.
Suppose that $\psi$ and $\psi'$ satisfy the conditions of the theorem.
Then $\psi\circ\varphi,\psi'\circ\varphi\colon H_2G\to A$ are coefficient
maps that map $\u_G$ to $v$; thus, these maps are equal, by the universality
property of $\u_G$. Since $\varphi$ is epimorphic, we deduce $\psi=\psi'$.\qed
\smallskip

\bigskip

\centerline{\bf III. $SL(2,K)$}
\bigskip
In this part we specialize our considerations to the discrete group $SL(2,K)$,
where $K$ is an infinite field.
This group is perfect, so that the results of Section 6 apply.
Our main result is that the (reduced) Witt class is the universal
equicommutative class for this group (Theorem 10.1). 
The proof relies on several known results which we review carefully.
The Schur multiplier of $SL(2,K)$, denoted $\fund{K}$ henceforth to honour
Calvin Moore, is classically described by generators and relations;
we recall this description in
Section 7. The generators are ``symbols'' $\{a,b\}$, $a,b\in\dotK$.
(We use $\dotK$ to denote $K\setminus\{0\}$.) In the quotient $\Eq(SL(2,K))$
of $\fund{K}$ the symbols become symmetric: $\{a,b\}=\{b,a\}$. We are thus led
to consider the group $\fund{K}/{\rm sym}$, defined by adjoining to the classical presentation of
$\fund{K}$ all the symbol symmetry relations $\{a,b\}=\{b,a\}$, $a,b\in\dotK$.
This group is
a natural mid-step
in the quotient sequence
$$\fund{K}\to\fund{K}/{\rm sym}\to\Eq(SL(2,K)).$$
In Section 8 we show that in fact $\fund{K}/{\rm sym}\simeq I^2(K)$
(here $I^2(K)$ is the square of the fundamental ideal $I(K)$ of the Witt
ring $W(K)$, cf.~[EKM, Chapter I]).
A cocycle $b$ representing the universal class $\u_{SL(2,K)}$ was given explicitly
(though slightly erroneously) by Moore;
we recall the correct description in Section 9.
There we also present the results of Nekov\'a$\check{\rm r}$,
and Kramer and Tent,
proving that the image of the Moore cocycle $b$ in $I^2(K)$
is cohomologous to the (reduced) Witt cocycle.
In Section 10 we use this compatibility to show that $I^2(K)\simeq\Eq(SL(2,K))$
and that the reduced Witt class is (equivalent to) the universal equicommutative class. 

\medskip
\bf 7. Schur multiplier of $SL(2,K)$.\rm
\medskip
\def\a{7}
\m=1
\n=1

In this section we recall the standard description of the universal central
extension and of the Schur multiplier of $SL(2,K)$. (As always, we assume that
$K$ is an infinite field.) The classical references are
[Moore, Sections 8,9], [Mats], [St, \S 7].

The universal central extension of $SL(2,K)$ is called the Steinberg group
and is denoted $St(2,K)$. It is generated by two families of symbols:
$x_{12}(t), t\in K$; $x_{21}(t), t\in K$. For $t\in\dotK$ one defines
an auxiliary element $w_{ij}(t)=x_{ij}(t)x_{ji}(-t^{-1})x_{ij}(t)$; then
the relations defining $St(2,K)$ are:
$$\eqalign{
x_{ij}(t)x_{ij}(s)&=x_{ij}(t+s)\qquad(t,s\in K),\cr
w_{ij}(t)x_{ij}(r)w_{ij}(t)^{-1}&=x_{ji}(-t^{-2}r)\qquad(t\in\dotK,r\in K),}
\leqno(\a.\f)$$
where $\{i,j\}=\{1,2\}$. Other noteworthy elements of $St(2,K)$ are
$h_{ij}(t)=w_{ij}(t)w_{ij}(-1)$.

The projection $\pi\colon St(2,K)\to SL(2,K)$ is defined by
$$\pi\colon\quad x_{12}(t)\mapsto\pmatrix{1&t\cr0&1},
\quad x_{21}(t)\mapsto\pmatrix{1&0\cr t&1}.\leqno(\a.\f)$$
Then one easily checks that
$$\pi\colon\quad w_{12}(t)\mapsto\pmatrix{0&t\cr-t^{-1}&0},
\quad h_{12}(t)\mapsto\pmatrix{t&0\cr 0&t^{-1}}.\leqno(\a.\f)$$

The kernel of $\pi$, i.e.~the Schur multiplier of $SL(2,K)$
(denoted $H_2(SL(2,K))$, $\fund{K}$ or $KSp_2(K)$ in various sources)
is generated by elements
$$\{s,t\}:=h_{12}(s)h_{12}(t)h_{12}(st)^{-1}\quad (s,t\in\dotK).\leqno(\a.\f) $$
With this generating set it is described by an explicit family of relations:
$$\eqalign{
&\{st,r\}\{s,t\}=\{s,tr\}\{t,r\},\quad \{1,s\}=\{s,1\}=1;\cr
&\{s,t\}=\{t^{-1},s\};\cr
&\{s,t\}=\{s,-st\};\cr
&\{s,t\}=\{s,(1-s)t\} \quad \hbox{\rm if $s\ne1$}.}\leqno(\a.\f)
$$
\noindent (For this presentation see [Moore, Theorem 9.2] or
[St, \S 7, Theorem 12]. Even though the group is abelian, the
convention is multiplicative.)

Notice that $h_{12}(s)$ and $h_{12}(t)$ are lifts to $St(2,K)$ of commuting elements
$\pmatrix{s&0\cr0&s^{-1}}$, $\pmatrix{t&0\cr0&t^{-1}}$. Therefore, the commutator
$[h_{12}(s),h_{12}(t)]=\{s,t\}\{t,s\}^{-1}$ belongs
to $\com({SL(2,K)})$---the kernel
of the quotient map $\fund{K}\to \Eq(SL(2,K))$. Imposing in $\fund{K}$
the extra ``symbol symmetry'' relations $\{s,t\}=\{t,s\}$ we obtain the group
$\fund{K}/{\rm sym}$, ``the symmetrized Schur multiplier''---an intermediate
step in passing from $\fund{K}$ to $\Eq(SL(2,K))$.
We will prove that this group is in fact equal to $\Eq(SL(2,K))$.
For this we use quadratic form theory.

\medskip
\bf 8. Symmetrized Schur multiplier and quadratic forms.\rm
\medskip
\def\a{8}
\m=1
\n=1

The fundamental ideal $I(K)$ of $W(K)$
is generated by non-degenerate \quadratic forms
on even-dimensional spaces. Another suitable collection of generators consists
of the forms $\pf(a)=\langle1,-a\rangle$, $a\in\dotK$.
The ideal $I^2(K)$ is the square of $I(K)$; it is generated by
Pfister forms $\pf(a,b)$ (for $a,b\in \dotK$), where
$\pf(a,b)=\pf(a)\otimes\pf(b)=\langle 1,-a\rangle\otimes\langle 1,-b\rangle=
[1]-[a]-[b]+[ab]$ (the last equality valid in $W(K)$).
More on these generating sets (in particular, the relations)
can be found in [EKM, I.4].

Let us state (in our current notation) [Su, Corollary 6.4]:
there exists a natural homomorphism $\Phi\colon\fund{K}\to I^2(K)$,
sending $\{a,b\}$ to $\pf(a,b)$; the kernel of $\Phi$ is generated
by the elements $\{a^2,b\}$. 
We will also need 
[Moore, Lemma 3.2]: in $\fund{K}$
$$\{a,b\}\{b,a\}^{-1}=\{a^2,b\}.\leqno(\a.\f)$$
Putting these two facts together we get the following.
\smallskip
\bf Proposition \a.\t.\sl\par
The homomorphism $$\Phi\colon\fund{K}\ni\{a,b\}\mapsto\pf(a,b)\in I^2(K)$$
induces an isomorphism $$\phi\colon\fund{K}/{\rm sym}\to I^2(K).$$
\rm
\smallskip
\bf Remark \a.\t.\rm\par\nobreak
It is known that, in general, the map $\Phi$ is not an isomorphism.
This means that for some field $K$ and some $a,b\in\dotK$
we have $\{a,b\}\ne\{b,a\}$ in $\fund{K}$.
It follows that for that field $K$ the Moore class is not equicommutative.
In topological terms, we may consider the $SL(2,K)$-bundle $\xi_{a,b}$ over the torus $T^2$
with monodromies
$\pmatrix{a&0\cr0&a^{-1}}$, $\pmatrix{b&0\cr0&b^{-1}}$.
Then the inequality $\{a,b\}\ne\{b,a\}$ implies that
the Moore class of $\xi_{a,b}$ is non-trivial.
\smallskip
\bf Remark \a.\t.\rm\par\nobreak
[Su, Corollary 6.4], crucial for the proof of Proposition \a.1,
follows from [Su, Lemma 6.3]. In the statement of that lemma the second
of the defining relations (2) is misprinted; the correct version is
$$\pf(a)+\pf(b)=\pf(a+b)+\pf({(a+b)ab}),\quad a+b\ne0.\leqno(\a.\f)$$
In the proof of part (3) of the lemma Suslin writes:
``It is trivial to check [\dots]
that the elements [\dots] satisfy relations (2)''.
To check Relation (\a.2) we needed the following
calculation in $\fund{K}/{\rm sym}$
(inspired by a 
calculation of Rost, cf.~[GSz, Lemma 7.6.8]):
$$
\eqalign{\{a,b\}&=\{b,a\}=\{b,-ba\}=\{-ab,b\}=\{-ab^{-1},b\}\cr&=
\{-ab^{-1},(1+ab^{-1})b\}=\{-ab,a+b\}=\{a+b,-ab\}=\{a+b,ab(a+b)\}.
}\leqno(\a.\f)$$
Each step uses symbol symmetry, applies one of the defining relations (7.5),
or multiplies a symbol entry by a square. The latter operation
is equivalent to multiplication by a 
symbol of the form
$\{z^2,y\}$,
as asserted in [Moore, Appendix, (7)]; by (\a.1), the symbol
$\{z^2,y\}$ is trivial in $\fund{K}/{\rm sym}$.

\medskip
\bf 9. Comparison of the Moore and Witt cocycles.\rm
\medskip
\def\a{9}
\m=1
\n=1

In this section we recall the explicit form of a cocycle $b$
representing the universal class $\u_{SL(2,K)}\in H^2(SL(2,K),\fund{K})$
as given in [Moore, 9.1-4], with later corrections (cf.~[Kr-T, 9.1]).
We also describe the image of $b$ under the map $\Phi\colon\fund{K}\to I^2(K)$.
Kramer and Tent show that 
this image, an $I^2(K)$-valued cocycle on $SL(2,K)$, is cohomologous to the Witt cocycle.
In the next section we will use this fact to show that the
$I^2(K)$-valued Witt class is the universal equicommutative class
for $SL(2,K)$.

Kramer and Tent do their calculation in the generic case, and argue
that this is enough to claim cocycle equality. We present the details in
all cases as this allows us to give an explicit formula for
a universal equicommutative cocycle.

Every element of $SL(2,K)$ is uniquely represented in one of the forms:
$$g_1(u,t)=x(u)h(t)=\pmatrix{1&u\cr0&1}\pmatrix{t&0\cr0&t^{-1}};\leqno(\a.\f)$$
$$g_2(u,t,v)=x(u)w(t)x(v)=\pmatrix{1&u\cr0&1}\pmatrix{0&t\cr -t^{-1}&0}
\pmatrix{1&v\cr0&1}.\leqno(\a.\f)$$
This leads to the following definition of a lift $SL(2,K)\to St(2,K)$:
$$\ov{g_1(u,t)}=\ov{x(u)h(t)}:=x_{12}(u)h_{12}(t);\leqno(\a.\f)$$
$$\ov{g_2(u,t,v)}=\ov{x(u)w(t)x(v)}:=x_{12}(u)w_{12}(t)x_{12}(v)
\leqno(\a.\f)$$
The corresponding cocycle $b$ was calculated by Moore, with later correction
by Schwarze (cf.~[Kr-T, 9.1]). We present the formulae for the cocycle $b$,
and for its image under $\Phi$ in $W(K)$.

\item{1)}
{ $\displaystyle{b(g_2(u,t,v),g_2(u',t',v'))=
\cases{
\{-w't^{-1}t'^{-1},-tt'^{-1}\}\{-t,-t'\}^{-1}& if $w':=-(v+u')\ne0$,\cr
\{-t,-t'\}^{-1}& if $w'=0$.
}}$

These are mapped by $\Phi$ to
$$\pf({-w't^{-1}t'^{-1},-tt'^{-1}})-\pf(-t,-t')=
\pf(-w'tt',-tt')-\pf(-t,-t')=[w']-[t]-[t']+[tt'w'].$$ 

and to
$$-\pf(-t,-t')=-[1]-[t]-[t']-[tt'].$$
}

\item{2)}
{ $b(g_2(u,t,v),g_1(u',t'))=\{t,t'^{-1}\}$.
This is mapped by $\Phi$ to $\pf(t,t'^{-1})=\pf(t,t')$.
}
\item{3)}
{ $b(g_1(u,t),g_2(u',t',v'))=\{t,t'\}$,
mapped to $\pf(t,t')$.
}
\item{4)}
{ $b(g_1(u,t),g_1(u',t'))=\{t,t'\}$,
mapped to $\pf(t,t')$.
}

We summarize:
$$(\Phi_*b)(g,h)=
\cases{
[w']-[t]-[t']+[tt'w']&
if $g=g_2(u,t,v)$, $h=g_2(u',t',v')$, $w'=-(v+u')\ne0$,\cr  
-[1]-[t]-[t']-[tt']&
if $g=g_2(u,t,v)$, $h=g_2(u',t',v')$, $w'=-(v+u')=0$,\cr
\pf(t,t')&
if $g=g_1(u,t)$, $h=g_2(u',t',v')$,\cr
\pf(t,t')&
if $g=g_2(u,t,v)$, $h=g_1(u',t')$,\cr
\pf(t,t')&
if $g=g_1(u,t)$, $h=g_1(u',t')$.
}\leqno(\a.\f)
$$

On the other hand, we have the Witt class $\w\in H^2(SL(2,K),W(K))$,
given by the cocycle $w$ defined by (1.4):
$$w(g,h)=[-|e,ge|\cdot|e,ghe|\cdot|e,he|],\leqno(\a.\f)$$
where $e={1\choose0}$. We now express this cocycle in the parametrization
of $SL(2,K)$ used by Moore.
\smallskip
\bf Lemma \a.\t.\sl\par
$$w(g,h)=\cases{
[w']& if $g=g_2(u,t,v)$, $h=g_2(u',t',v')$, $w'=-(v+u')\ne0$,\cr  
0& in all other cases.}
\leqno(\a.\f)$$
\rm\par
Proof.
Notice that
$g_1(u,t)e=x(u)h(t)e={t\choose0}$, so that $|e,g_1(u,t)e|=0$.
Therefore, 
the value of $w(g,h)$ is zero if any of the arguments $g,h$ is of the form
$g_1(u,t)$.
Even so,
for later use we need the following easily checked formulae:
$$
g_2(u,t,v)g_1(u',t')e={*\choose -t^{-1}t'},\quad
g_1(u,t)g_2(u',t',v')e={*\choose -t^{-1}t'^{-1}},\quad
g_1(u,t)g_1(u',t')e={tt'\choose0}.\leqno(\a.\f)$$

Let us turn to case 1:
$$g_2(u,t,v)e=x(u)w(t)x(v)e=x(u)w(t)e=x(u){0\choose -t^{-1}}
={-t^{-1}u\choose-t^{-1}},\leqno(\a.\f)$$
so that $|e,g_2(u,t,v)e|=-t^{-1}$. Furthermore (setting $w'=-v-u'$),
$$\eqalign{
g_2(u,t,v)g_2(u',t',v')e&=
x(u)w(t)x(v){-t'^{-1}u'\choose-t'^{-1}}=
x(u)\pmatrix{0&t\cr-t^{-1}&0}{t'^{-1}w'\choose -t'^{-1}}\cr
&=\pmatrix{1&u\cr0&1}{-tt'^{-1}\choose  -t^{-1}t'^{-1}w'}=
{-tt'^{-1} -t^{-1}t'^{-1}uw'\choose  -t^{-1}t'^{-1}w'}
},$$
so that
$$|e,g_2(u,t,v)g_2(u',t',v')e|=
-t^{-1}t'^{-1}w'.
\leqno(\a.\f)$$
Therefore
$$w(g_2(u,t,v),g_2(u',t',v'))=\cases{
[-(-t^{-1})(-t'^{-1})(-t^{-1}t'^{-1}w')]=[w']& if $w'\ne0$,\cr
0 & if $w'=0$.}\leqno(\a.\f)
$$
\hfill\qed
\smallskip

Nekov\'a$\check{\rm r}$ [Ne, \S 2] noticed that $w$ is cohomologous
to an $I^2(K)$-valued cocycle; it is enough to add the coboundary
of the following cochain:
$$n(g)=\cases{
[|e,ge|] & if $|e,ge|\ne0$,\cr
[1]-[t]& if $ge=te$ for some $t\in\dotK$. 
}\leqno(\a.\f)$$
In the parametrization used by Moore:
$n(g_1(u,t))=[1]-[t]$, $n(g_2(u,t,v))=[-t^{-1}]=-[t]$,
$$\displaylines{
n(g_2(u,t,v)g_2(u',t',v'))=\cases{
[-t^{-1}t'^{-1}w']=-[tt'w']& if $w'\ne0$,\cr
[1]-[-tt'^{-1}]=[1]+[tt']& if $w'=0$;
}\cr
(\a.\f)\hfill n(g_2(u,t,v)g_1(u',t'))=[-t^{-1}t']=-[tt'],\qquad
n(g_1(u,t)g_2(u',t',v'))=[-t^{-1}t'^{-1}]=-[tt'],\cr
n(g_1(u,t)g_1(u',t'))=[1]-[tt'].}$$
We are ready to calculate $\delta n$ and see that $\Phi_*b=w+\delta n$.
Using the formula
$(\delta n)(g,h)=n(g)-n(gh)+n(h)$ we get:
$$
\leqalignno{
(\delta n)(g_2(u,t,v),g_2(u',t',v'))
&=\cases{
-[t]+[tt'w']-[t'] & if $w'=-(v+u')\ne0$,\cr
-[t]-([1]+[tt'])-[t'] & if $w'=0$.
}&1)\cr
(\delta n)(g_2(u,t,v),g_1(u',t'))&
=-[t]-(-[tt'])+[1]-[t']=\pf(t,t')&2)\cr
(\delta n)(g_1(u,t),g_2(u',t',v'))&
=[1]-[t]-(-[tt'])-[t']=\pf(t,t')&3)\cr
(\delta n)(g_1(u,t),g_1(u',t'))&
=[1]-[t]-([1]-[tt'])+[1]-[t']=\pf(t,t')&4)}$$
Comparing these four formulae, (\a.5) and (\a.7) we obtain the following
proposition.

\smallskip
\bf Proposition  \a.\t.\rm ([Kr-T, 9.2]) \sl\par
The Witt cocycle $w$ and the cocycle $\Phi_*b$ are cohomologous
(in the complex of  $SL(2,K)$-cochains with coefficients in $W(K)$).
\rm\par
\bigskip
The phrasing of Proposition \a.2 is slightly awkward, because
$\Phi_*b$ has coefficients in $I^2(K)$, while $w$ has coefficients in $W(K)$.
Fortunately, it is not hard to check that the relevant cohomology groups embed:
\smallskip
\bf Proposition \a.\t.\sl\par
The inclusion $\iota\colon I^2(K)\to W(K)$ induces a monomorphism
$$\iota_*\colon H^2(SL(2,K),I^2(K))\to H^2(SL(2,K),W(K)).\leqno(\a.\f)$$
\rm\par
Proof. Let $Q=W(K)/I^2(K)$.
Consider the short exact sequence of coefficient groups:
$$0\to I^2(K)\to W(K)\to Q\to0\leqno(\a.\f)$$
and the associated long exact sequence
$$\ldots\to H^1(SL(2,K),Q)\to H^2(SL(2,K),I^2(K))\to H^2(SL(2,K),W(K))
\to\ldots\leqno(\a.\f)$$
We have
$$H^1(SL(2,K),Q)=\Hom(SL(2,K),Q)=0,\leqno(\a.\f)$$
because $SL(2,K)$ is perfect, hence has no nontrivial homomorphisms
to abelian groups.\qed
\smallskip
To summarize: $\iota_*([\Phi_*b])=[w]$---or, in terms of cohomology classes,
$\iota_*(\Phi_*(\u_{SL(2,K)}))=\w$.
\smallskip
\bf Definition \a.\t.\rm\par
The \it reduced \rm (or \it $I^2(K)$-valued \rm\hskip-.1cm) Witt class $\w^I\in H^2(SL(2,K),I^2(K))$
is defined by $\w^I=\iota_*^{-1}(\w)$; it is equal to $\Phi_*(\u_{SL(2,K)})$
and represented by the cocycle $\Phi_*b$,
explicitly given by (\a.5).
\smallskip
For practical purposes, one can ignore the difference between the classes $\w$ and $\w^I$,
mainly because of the following corollary of Proposition \a.2.
\smallskip
\bf Corollary \a.\t.\sl\par
For any $SL(2,K)$-bundle $P$ over a closed oriented surface $\Sigma$
we have
$$\langle \w(P),[\Sigma]\rangle=
\langle {\w}^I(P),[\Sigma]\rangle\in I^2(K).
\leqno(\a.\f)$$
\smallskip
\bf Proposition \a.\t.\sl\par
The class $\w^I$ is equicommutative.
\rm\par
Proof. A cocycle $c$ representing $\w^I$ treated as a $W(K)$-valued cocycle
(via the embedding $I^2(K)\to W(K)$) is cohomologous to the (standard)
Witt cocycle $w$. The latter is equicommutative, hence, by Proposition 5.3,
so is $c$. \qed

\medskip
\bf 10. The Witt class is universal equicommutative.\rm
\medskip
\def\a{10}
\m=1
\n=1

We prove what is in the title of this section.
\smallskip
\bf Theorem \a.\f. \rm(Theorem A)\sl\par
Let $K$ be an infinite field.
The group $\Eq(SL(2,K))$ is isomorphic to $I^2(K)$, and
the Witt class $\w^I\in H^2(SL(2,K),I^2(K))$ is the universal equicommutative
class.
\rm\par
Proof.
Consider the following diagram of coefficient groups.
$$
\matrix{
\fund{K}
&\mathop{\longrightarrow}\limits^{q_1}
&\fund{K}/{\rm sym}
&\mathop{\longrightarrow}\limits^{q_2}
&\Eq(SL(2,K))\cr
&{d_1}\searrow\phantom{D}&
\Big{\downarrow}{\phi}
&\phantom{D}\swarrow{d_2}&
\cr
&
&
I^2(K)
&
&
}\leqno(\a.\f)
$$
The diagonal arrows are the unique maps deduced from the universal properties
of $\u_{SL(2,K)}$ and $\ww_{SL(2,K)}$, applied to the Witt class $\w^I$.
Uniqueness of the universal map implies commutativity of the diagram.
Namely, the left triangle commutes because $\phi q_1=\Phi$ maps
$\u_{SL(2,K)}$ to $\w^I$ (Proposition 9.2), hence is equal to the diagonal map $d_1$.
Similarly, we have $\w^I=(d_2)_*\ww_{SL(2,K)}=(d_2)_*(q_2q_1)_*\u_{SL(2,K)}$, hence
$d_2q_2q_1=d_1=\phi q_1$; but $q_1$ is surjective, therefore we deduce
$d_2q_2=\phi$, i.e.~the commutativity of the right triangle.
Now $d_2$ 
is an isomorphism, because
$\phi$ is an isomorphism (Proposition 8.1) and $q_2$ is surjective.
The isomorphism $d_2$ maps the universal equicommutative
class $\ww_{SL(2,K)}$ to $\w^I$. 
\qed
\smallskip
Another corollary of the proof is that $q_2$ is an isomorphism:
the quotient of $\fund{K}$ by the symbol symmetry relations
is already equal to $\Eq(SL(2,K))$.

\bigskip
\centerline{\bf IV. Witt range.}
\bigskip
What are the possible values of the Witt class for
$SL(2,\Q)$-bundles over surfaces? We know that these values reside in
$I^2\Q$, which is a direct sum of $4\Z$ (the real, signature part) and an
infinite direct sum of $\Z/2$ ($p$-adic parts, one per odd prime). More details of this description
are given in Section 12. While the real part of the Witt class
is related to the Euler class, hence non-trivial (cf.~[\TccI, Section 13]),
the $p$-adic parts are more mysterious---perhaps trivial? Using (1.4)
we did some computer calculations in FriCAS that indicated non-triviality
of the $p$-adic parts. In Section 12 we give a complete description
of the range of the Witt class over $\Q$, proving that the Milnor--Wood inequality
(restricting the real part in terms of the genus of the base surface) is the
only restriction---the $p$-adic parts can be arbitrary. The challenging part of this
result is the construction of sufficiently many non-trivial bundles.
This is done in Section 11 in much greater generality, for arbitrary infinite
fields $K$. 
We analyse $SL(2,K)$-bundles
over simple surfaces (pair of pants; genus 1 surface with one boundary component)
using some results on the Markov equation (quoted from [GMS]). 
Then we use the gluing results from Section 4. The formula (1.4) is used
throughout to control the Witt class. Our final result, Theorem 11.6,
gives a large subset of Witt classes in 
$I^2(K)$.
This subset is quite close to the one
defined by the boundedness restriction for the Witt class; the difference
is discussed in Section 13.

An early paper where 
many $SL(2,\Q)$-bundles were constructed by gluing
is [Takeuchi].
\medskip
\bf 11. Markov surfaces and representations.\rm
\medskip
\def\a{11}
\m=1
\n=1

In this section  $K$ will be
an arbitrary infinite field. 
Let $\w$ be the Witt class. We will use the representing cocycle
$$w(X,Y)=[-|e,Xe|\cdot|e,XYe|\cdot|e,Ye|]=[-X_{21}\cdot(XY)_{21}\cdot Y_{21}],
\leqno(\a.\f)$$
where $e={1\choose0}$,  $X_{21}$ denotes the $21$-entry of the matrix $X$, and $[0]$
is interpreted as 0. Additive convention will be used for the cocycle.

\smallskip
\bf Remark \a.\t.\rm\par
If either $X$, $Y$ or $XY$ is diagonal (or even upper-triangular), then $w(X,Y)=0$.
\smallskip
We will use notions and notation discussed at the beginning of Section 4,
in particular the notion of standard framing,
and the notion of ``relative class''
$\ov{w}$ (cf.~Definition 4.1).
\smallskip
\bf Lemma \a.\t.\sl\par
Let $K$ be an infinite field. 
Then, for any $\alpha,\beta\in\dotK$ there exists a flat
$\G$-bundle $\xi$ with a standard framing, over the oriented
genus 1 compact surface with one boundary component, such that
$\ov{w}(\xi)=[\alpha]+[\beta]\in W(K)$.
Moreover, for every $z\in-\alpha\beta\dotKK$ (with finitely many exceptions)
the bundle $\xi$ may be chosen so that its boundary monodromy
is diagonal with eigenvalues $z,z^{-1}$.
\rm\smallskip
Proof. To construct the bundle as in the lemma, with boundary monodromy $Z$,
we need to find $X,Y\in \G$ such that $[X,Y]=Z$.
We heavily rely on the classical description of the space of solutions
of the commutator equation $[X,Y]=Z$; we use the version from [GMS], though
some results go back as early as to Fricke. To start, if $[X,Y]=Z$, then
the Fricke identity says that
the scalars $x_1=\tr{X}$, $x_2=\tr{Y}$, $x_3=\tr{XY}$ and $m=\tr{Z}+2$
satisfy the 
Markov equation
$$x_1^2+x_2^2+x_3^2-x_1x_2x_3=m.\leqno(M_m)$$
Conversely, if
\item{(a)} $\tr{Z}\ne\pm2$;
\item{(b)} $(x_1,x_2,x_3)$ is a solution of $(M_m)$;
\item{(c)} $m-x_2^2\ne0$;
\item{(d)} {$Y\in \G$ satisfies $\tr{Y}=\tr{ZY}=x_2$;}

\noindent then there exists a unique $X\in \G$ that satisfies $\tr{X}=x_1$,
$\tr{XY}=x_3$ and $[X,Y]=Z$ (cf.~[GMS, Lemma 3.5]).
Explicitly, this $X$ is given by
$$(m-x_2^2)X=-x_3ZY+x_1Z+(x_3-x_1x_2)Y^{-1}+x_1I.\leqno(\a.\f)$$
We will work out the case $Z=\pmatrix{z&0\cr0&z^{-1}}$, $z\ne\pm1$
(then (a) is fulfilled). The matrices $Y$ satisfying (d) can be found explicitly.
If $Y=\pmatrix{a&*\cr*&d}$, then $\tr{Y}=\tr{ZY}=x_2$ are equivalent to
a linear system of equations on $a$, $d$ with unique solution
$$a={x_2\over 1+z},\qquad d={zx_2\over1+z}.\leqno(\a.\f)$$
The condition $\det{Y}=1$ gives the following form of $Y$:
$$Y={1\over1+z}
\pmatrix{x_2&c^{-1}(zx_2^2-(1+z)^2)\cr c&zx_2}.\leqno(\a.\f)$$
The solution depends on a unique parameter $c\in\dotK$.
(We have $c\ne0$, since otherwise $1=\det{Y}$ would imply $1=zx_2^2(1+z)^{-2}=x_2^2m^{-1}$, contradicting
(c).)
Condition (c) ensures that $m-x_2^2\ne0$; therefore we may
plug (\a.4) into (\a.2) and determine $X$:
$$(m-x_2^2)X={1\over 1+z}
\pmatrix{x_1(1+z)^2-x_1x_2^2z&*\cr(x_1x_2-(1+z^{-1})x_3)c&x_1(1+z)(1+z^{-1})-x_1x_2^2}.\leqno(\a.\f)$$
Further direct computations show that
$$(m-x_2^2)(XY)_{21}={x_1(1+z)-x_2x_3\over z(1+z)}c,\qquad
(m-x_2^2)(YX)_{21}={x_1(1+z)-x_2x_3\over 1+z}c.\leqno(\a.\f)$$
Let $\xi_c$ be the bundle defined by these $(X,Y)$.
Finally, using Lemma 4.8 and Remark \a.1 we get
$$
\eqalign{
\ov{w}(\xi_c)
&=
w(X,Y)-w(Y,X)-w(Z,YX)\cr
&=[-X_{21}\cdot(XY)_{21}\cdot Y_{21}]-[-X_{21}\cdot(YX)_{21}\cdot Y_{21}]\cr
&=
[-{c\over z(1+z)}(x_1x_2-(1+z^{-1})x_3)(x_1(1+z)-x_2x_3)]\cr
&\phantom{=}-[-{c\over 1+z}(x_1x_2-(1+z^{-1})x_3)(x_1(1+z)-x_2x_3)].
}
\leqno(\a.\f)
$$
If our solution $(x_1,x_2,x_3)$ satisfies two further conditions

\item{(e)} $x_1x_2-(1+z^{-1})x_3\ne0$, $x_1(1+z)-x_2x_3\ne0$;

\noindent then we modify $c$ to a new parameter
$$C:=-{c\over z(1+z)}(x_1x_2-(1+z^{-1})x_3)(x_1(1+z)-x_2x_3)\leqno(\a.\f)$$
($C$ runs through $\dotK$ as $c$ does), and then we get
$$\ov{w}(\xi_c)=[C]-[Cz]=[C]+[-Cz].\leqno(\a.\f)$$
To summarize, to realize $[\alpha]+[\beta]\in W(K)$ (for given $\alpha,\beta\in\dotK$)
as $\ov{w}(\xi_c)$ we may choose $z=-\alpha\beta\lambda^2$, $C=\alpha$ 
(for some $\lambda\in\dotK$) and apply the above 
construction---provided that we find a solution
$(x_1,x_2,x_3)$ of $M_m$ that satisfies (c) and (e). In [GMS, proof
of Proposition 3.6] the following solution of $M_{t+2}$ is given:
$$x_2=\zeta+\zeta^{-1},\qquad x_3={-x_2^2+t+1\over \zeta-\zeta^{-1}},
\qquad x_1=1+\zeta x_3.\leqno(S)$$
(We have corrected a sign mistake in $x_3$; $\zeta\in\dotK\setminus\{\pm1\}$
is arbitrary.) Putting $t=z+z^{-1}$ we get:
$$x_1={\zeta z-\zeta(\zeta^2+2\zeta^{-2})+\zeta z^{-1}\over \zeta-\zeta^{-1}} ,\qquad
x_2={\zeta^2-\zeta^{-2}\over \zeta-\zeta^{-1}},\qquad
x_3={z-(\zeta^2+1+\zeta^{-2})+z^{-1}\over \zeta-\zeta^{-1}}
\leqno(S')$$

Let us now discuss the conditions (c) and (e) for the solution $(S')$. 
For a fixed $\zeta$,
these are three non-trivial polynomial inequalities on $z$, 
hence they hold except for some finite set $E$ of values of $z$.
(Since the general procedure described above requires $z\ne\pm1$, we
include these two values in $E$ as well.)
Clearly, for all but finitely many values of $\lambda$ we have $z\not\in E$;
this proves the lemma.
\qed
\smallskip
\bf Lemma \a.\t.\sl\par
Let $\lambda_1,\lambda_2,\lambda_3\in\dotK\setminus\{\pm1\}$. Let 
$$\Lambda_i=
\pmatrix
{
\lambda_i & 0 \cr 
0 & \lambda^{-1}_i
},
\qquad
t_i=\tr{\Lambda_i}=\lambda_i+\lambda^{-1}_i.\leqno(\a.\f)$$
Suppose that the triple $(t_1,t_2,t_3)$ does not satisfy Markov's equation $M_4$.
Then for every $c\in\dotK$ there exist a unique pair of matrices
$(L,M)$ in $\G$ that satisfies the following conditions:
\item{(a)} $L={}^A\Lambda_1$, $M={}^B\Lambda_2$ for some $A,B\in \G$;
\item{(b)} $LM=\Lambda_3$;
\item{(c)} {$L_{21}=c$.}

\noindent Moreover, regardless of the choice of $A$, $B$ 
in condition (a), we then have 
$$w({}^A\Lambda_1,A)+w({}^B\Lambda_2,B)=
[-c\lambda_1]+[c\lambda_2\lambda_3].\leqno(\a.\f)$$
\noindent (Recall that ${}^A\Lambda$ is our notation for the conjugate $A\Lambda A^{-1}$.)
\smallskip\rm
\bf Corollary \a.\t.\sl\par
Let $\xi$, $\xi'$ be flat $\G$-bundles with standard boundary framings and
diagonal boundary monodromies
$\pmatrix{\lambda_1&0\cr0&\lambda^{-1}_1},
\pmatrix{\lambda_2&0\cr0&\lambda^{-1}_2}
\ne\pm I$. Let 
$\alpha,\beta\in\dotK$. Then, for every
$z\in-\alpha\beta\lambda_1\lambda_2\dotKK$ (with finitely many exceptions)
there exist matrices $A,B\in \G$ such that
${}^A\xi\vee{}^B\xi'$ has diagonal boundary monodromy $
\pmatrix{z&0\cr0&z^{-1}}\ne\pm I$
and with standard boundary framing satisfies
$$\ov{w}({}^A\xi\vee{}^B\xi')=[\alpha]+[\beta]+\ov{w}(\xi)+\ov{w}(\xi').
\leqno(\a.\f)$$

\rm\par
Proof (of the Corollary).
We use Lemma \a.3 with $\Lambda_1$, $\Lambda_2$ equal to the boundary
monodromies of $\xi$, $\xi'$. We choose 
$\lambda_3=z=-\alpha\beta\lambda_1\lambda_2\lambda^2$ and $c=\alpha\lambda_1$
for some
$\lambda\in\dotK$. (There is a finite set of values that
$\lambda_3$ has to avoid:
$\pm1$ and the values for which $(t_1,t_2,t_3)$ would
satisfy $M_4$; we avoid them for all but finitely many choices of
$\lambda\in\dotK$.)
To prove (\a.12) we use Lemmas 4.6 and 4.7 (that describe the change of $\ov{w}$
under twists
and $\vee$), the remark about vanishing of $w$ for diagonal matrices, and (\a.11):
$$\eqalign{
\ov{w}({}^A\xi\vee{}^B\xi')&=\ov{w}({}^A\xi)+\ov{w}({}^B\xi')+w({}^A\Lambda_1,{}^B\Lambda_2)\cr
&=\ov{w}(\xi)+w(A,\Lambda_1)-w({}^A\Lambda_1,A)+\ov{w}(\xi')+w(B,\Lambda_2)-w({}^B\Lambda_2,B)\cr
&=\ov{w}(\xi)+\ov{w}(\xi')-(w({}^A\Lambda_1,A)+w({}^B\Lambda_2,B))\cr
&=\ov{w}(\xi)+\ov{w}(\xi')-([-c\lambda_1]+[c\lambda_2\lambda_3])\cr
&=\ov{w}(\xi)+\ov{w}(\xi')+[\alpha]+[\beta].}
\leqno(\a.\f)$$
\qed
\smallskip
Proof (of Lemma \a.3).
Suppose that $(L,M)$ is a pair satisfying (a), (b) and (c).
Let $$L=\pmatrix{a&d\cr c&t_1-a}.$$
We claim that $c\ne0$. Indeed, if $c$ were $0$ then
$L$ would be upper-triangular (with eigenvalues $\lambda^{\pm1}_1$),
and $M=L^{-1}\Lambda_3$ would also be upper-triangular
(with eigenvalues $\lambda^{\pm1}_2$). Direct calculation shows that then
$[L,M]$ would be upper-triangular and of trace $2$
(with both diagonal entries equal to $1$). The Fricke identity for
the traces $(t_1=\tr{L},t_2=\tr{M},t_3=\tr{\Lambda_3})$ would
show that this triple satisfies the Markov equation $(M_4)$---contradiction.

Thus, we have $L=\pmatrix{a&d\cr c&t_1-a}$ for some $c\in\dotK$.
It follows that 
$d=-(a^2-t_1a+1)/c$.
Next:
$$M=L^{-1}\Lambda_3=
\pmatrix{t_1-a&-d\cr -c&a}
\pmatrix{\lambda_3&0\cr 0&\lambda^{-1}_3}
=
\pmatrix{\lambda_3(t_1-a)&-\lambda^{-1}_3d\cr -\lambda_3c&\lambda^{-1}_3a}.
\leqno(\a.\f)
$$
The trace of $M$ is $t_2$:
$\lambda_3(t_1-a)+\lambda^{-1}_3a=t_2$, which holds for
$$a={\lambda_3t_1-t_2\over\lambda_3-\lambda^{-1}_3}.\leqno(\a.\f)$$
To finish the discussion of the existence and uniqueness question
we remark that any matrix with determinant 1 and trace $t_i$ is conjugate
(in $\G$) to $\Lambda_i$.

We now pass to the calculation of $w$. The equation ${}^A\Lambda_1=L$
implies that $A_1$, the left column of $A=(A_1,*)$, is a $\lambda_1$-eigenvector
of $L$. Therefore $$(LA)_{21}=(\lambda_1A_1,*)_{21}=\lambda_1A_{21},$$
so that
$$w(L,A)=[-L_{21}(LA)_{21}A_{21}]=[-c\cdot\lambda_1A_{21}\cdot A_{21}]=
[-c\lambda_1].
\leqno(\a.\f)
$$
(We know that $A_{21}\ne0$, for otherwise $A$ and $^A\Lambda_1=L$ would be upper-triangular,
contradicting $L_{21}=c\ne0$.)
Similarly,
$$w(M,B)=[-M_{21}(MB)_{21}B_{21}]=
[-(-\lambda_3c)\cdot(\lambda_2B_{21})\cdot B_{21}]=[c\lambda_2\lambda_3].
\leqno(\a.\f)
$$
\qed
\smallskip

One of the ways to look at the Witt ring $W(K)$ is the following.
Any element of $W(K)$ is represented by a unique (up to isomorphism)
anisotropic form over $K$. The dimension of this anisotropic representative
defines a norm $\|\cdot\|$ on $W(K)$ (cf.~[MH, 3.1.7, 3.1.8]).
On the other hand, any element $x\in W(K)$ can be represented---in many ways---as a finite
sum
$$x=\sum_{i\in I}n_i[a_i],\leqno(\a.\f)$$
where $n_i\in\Z$, $a_i\in\dotK$. The symbolic norm $\|x\|_s$ is the
minimum---over all such represen\-ta\-tions---of the expression $\sum_{i\in I}|n_i|$.
\smallskip
\bf Lemma \a.\t.\sl\par
For each $x\in W(K)$ we have $\|x\|=\|x\|_s$.
\smallskip\rm
Proof. If $x=\sum_{i\in I}n_i[a_i]$, then the form $x=\bigoplus_{i\in I}|n_i|\langle(\sgn{n_i})a_i\rangle$
represents $x$. This form has a largest anisotropic direct summand---the unique (up to isomorphism)
representative of $x$, of dimension $\|x\|$. It follows that
$$\|x\|\le\dim(\bigoplus_{i\in I}|n_i|\langle(\sgn{n_i})a_i\rangle)=\sum_{i\in I}|n_i|.
\leqno(\a.\f)$$
Therefore $\|x\|\le\|x\|_s$. On the other hand, a diagonalization of the anisotropic representative of $x$
expresses $x$ as a sum of $\|x\|$ symbols $[a_i]$, so that $\|x\|_s\le\|x\|$.\qed
\smallskip
The diameter of $W(K)$ with respect to the above norm is sometimes called
the $u$-invariant of $K$ (cf.~[Lam, XI.6]). However, it is too often infinite (e.g.~for $K=\Q$), hence
a refinement is widely used: $u(K)$ is the diameter of the set $W_t(K)$ of torsion elements in $W(K)$
(cf.~[EKM, Chapter VI] or [Lam, Definition XI.6.24]). It is classically known that $u(K)=4$ for
local (non--archimedean) and for global fields (cf.~[Lam, Examples XI.6.2, XI.6.29] or [EKM, Example 36.2]).

We use the following facts from [Lam, II.2]. The determinant of a quadratic form $q=\langle a_1,\ldots,a_n\rangle$
is defined as  $d(q):=\prod_{i=1}^n a_i\in\dotK/\dotKK$, and the discriminant as
$d_\pm(q)=(-1)^{n(n-1)/2}d(q)$. An even--dimensional form $q$ is in $I^2(K)$ if and only if $d_\pm(q)=1$.
\smallskip
\bf Theorem \a.\t. \rm(Theorem B)\sl\par
Let $K$ be an infinite field.
\item{(a)} The Witt class of any flat $SL(2,K)$-bundle over an oriented
closed surface of genus $g$ has norm $\le 4(g-1)+2$.
\item{(b)} 
The set of Witt classes of flat $SL(2,K)$-bundles over an oriented
closed surface of genus $g$ contains the set of elements
of $I^2(K)$ of norm $\le 4(g-1)$.
\smallskip\rm
Proof.

(a) This part is straightforward. The closed orientable genus $g$ surface $\Sigma_g$ has
a $\Delta$-complex structure with $4g-2$ triangles (cf.~the proof of Lemma 3.1).
The value of the Witt cocycle $w(\xi)$ (of any flat $SL(2,K)$-bundle $\xi$ over $\Sigma_g$)
on each of these triangles has norm $\le1$. This implies claim a).

(b) Let $q\in I^2(K)$ be an element of norm $\le 4g-4$. Then we can find
$\alpha_i,\beta_i,\gamma_j,\delta_j\in\dotK$ ($1\le i \le g$, $2\le j\le g-1$)
such that
$$q=\sum_i([\alpha_i]+[\beta_i])+\sum_j([\gamma_j]+[\delta_j]).\leqno(\a.\f)$$
If the norm of $q$ is smaller than $4g-4$, we add some extra trivial
terms $[1]+[-1]$ to obtain the above form;
since $q\in I^2(K)$ we know that $\dim{q}$ is even, and that the following
product formula holds:
$$1=d_\pm(q)=\prod_i(\alpha_i\beta_i)\prod_j(\gamma_j\delta_j).\leqno(\a.\f)$$
Now we find, using Lemma \a.2, bundles $\xi_i$ (for $i<g$),
(over genus 1 oriented surfaces with one boundary component),
with standard boundary framing,
with $\ov{w}(\xi_i)=[\alpha_i]+[\beta_i]$, 
with diagonal boundary monodromy with eigenvalues $z_i^{\pm1}$ that
satisfy $[z_i]=[-\alpha_i\beta_i]$.
We put $\zeta_1=\xi_1$, $u_1=z_1$, and then, using Cor.~\a.4, we
recursively define
$$\zeta_j:={}^{A_j}\zeta_{j-1}\vee{}^{B_j}\xi_j,\leqno(\a.\f)$$
with standard boundary framing, diagonal boundary monodromy
with eigenvalues $u_j^{\pm1}$ that satisfy
$$[u_j]=[-\gamma_j\delta_ju_{j-1}z_j]\leqno(\a.\f)$$
and
$$\ov{w}(\zeta_j)=[\gamma_j]+[\delta_j]+\ov{w}(\zeta_{j-1})+\ov{w}(\xi_j).
\leqno(\a.\f)$$
Induction then gives:
$$[u_j]=[-\prod_{i=1}^j(\alpha_i\beta_i)
\prod_{\ell=2}^j(\gamma_\ell\delta_\ell)],\leqno(\a.\f)$$
$$\ov{w}(\zeta_j)=\sum_{i=1}^j([\alpha_i]+[\beta_i])+
\sum_{\ell=2}^j([\gamma_j]+[\delta_j]).\leqno(\a.\f)$$
Now we construct $\xi_g$ just as the other $\xi_i$, but with
$\ov{w}(\xi_g)=[-\alpha_g]+[-\beta_g]$.
For $j=g-1$ we obtain, by the product formula, $[u_{g-1}]=[z_g]$.
It follows that $\zeta_{g-1}$ and $\xi_g$ can be constructed
with the same boundary monodromy. Then $\zeta_{g-1}\cup\xi_g$
is the desired bundle:
$$\w(\zeta_{g-1}\cup\xi_g)=\ov{w}(\zeta_{g-1})-\ov{w}(\xi_g)=q.\leqno(\a.\f)$$
\qed

\medskip
\bf 12. The range of the Witt class over $\Q$.\rm
\def\a{12}
\m=1
\n=1
\medskip
The goal of this section is to describe the range of the Witt class
for all representations of $\pi_1(\Sigma_g)$ (the fundamental group
of the genus $g$ orientable surface) in $SL(2,\Q)$.
\smallskip

We now recall a description of $I^2\Q$ from [Lam, VI.5.8], slightly modified using [MH, 4.2.5]: the sequence
$$0\to I^2\Q\to I^2\R\oplus\bigoplus_{p\ne\infty}I^2\Q_p\mathop{\rightarrow}\limits^{f}\Z/2\to0\leqno(\a.\f)$$
is exact. Here the middle terms are: $I^2\R\simeq 4\Z$, $I^2\Q_p\simeq \Z/2$. The first embedding is
defined by functorial maps associated to the completion embeddings $\Q\to\R$, $\Q\to\Q_p$, while
$f$ is the ``reciprocity law'' map:
$$f(4n,a_2,a_3,a_5,\ldots)=(n+\sum_{p\ne\infty}a_p) \mod 2.\leqno(\a.\f)$$
(Each $I^2\Q_p$ is isomorphic to $\Z/2$; each non-trivial $a_p$ is interpreted as $1$.)
The torsion part of $I^2\Q$ is isomorphically mapped onto the subgroup
$$\{(a_p)_{p\ne\infty}\in\bigoplus_{p\ne\infty}I^2\Q_p\mid \sum_{p\ne\infty}a_p\equiv 0\mod 2\}.\leqno(\a.\f)$$
We denote by $\sigma\colon W(\Q)\to W(\R)\simeq\Z$ the signature map, normalized by
$\sigma([1])=1$. Then the torsion elements of $W(\Q)$ are the ones of signature zero.
\smallskip
\bf Lemma \a.\t.\sl\par
\item{(a)} Elements of $I^2\Q$ of signature $\pm4 h$  (where $h\ge1$)
have norm $4h$. The set of such elements
can be described as
$$\{\sum_{i=1}^{4h}[a_i]\in W(\Q)\mid \pm a_i>0, \prod_ia_i\in\dotQQ\}.
\leqno(\a.\f)$$
\item{(b)} Non-trivial elements of $I^2\Q$ of signature $0$ have norm $4$. The set of such elements
can be described as
$$\{[a]+[b]+[c]+[d]\in W(\Q)\mid abcd\in\dotQQ, \hbox{\rm and exactly two of $a,b,c,d$ are positive}\}.\leqno(\a.\f)$$
\smallskip\rm
Proof.
Let us start with the proof of part (b). Let $x\in I^2\Q$, $x\ne0$, $\sigma(x)=0$. Since $u(\Q)=4$ we know
that $\|x\|\le4$. As the dimension function takes even values on $I\Q$ (by definition), hence also
on $I^2\Q$, an anisotropic representative of $x$ is 4- or 2-dimensional. If $x=\langle a,b\rangle$,
however, we get $1=d_\pm(\langle a,b\rangle)=-ab$; then $ab=-1$ in $\dotQ/\dotQQ$, and
$x=\langle a, -a\rangle=0$ in $W(\Q)$. Therefore, an anisotropic representative of $x$ is 4-dimensional.
Thus, we have $x=[a]+[b]+[c]+[d]$ for some $a,b,c,d\in\dotQ$ that satisfy
$abcd\in\dotQQ$ (this is equivalent to $x\in I^2\Q$),
and exactly two of $a,b,c,d$ are positive (equivalent to $\sigma(x)=0$).

For part (a) we use Meyer's theorem (cf.~[MH, Corollary II.3.2]): 
a quadratic $\Q$-form of dimension greater than 4 is $\Q$-isotropic if 
it is $\R$-isotropic. Let $x\in I^2\Q$, $\sigma(x)=4h>0$.
Let $q$ be an anisotropic
representative of $x$. If $\dim{q}>4h$, then $q$ is $\R$-isotropic,
hence, by Meyer's theorem,
also $\Q$-isotropic---contradiction.
The rest of the statement is seen as in part (b):
the condition $\prod_ia_i\in\dotQQ$, i.e.~$d_\pm=1$,
characterizes elements in $I^2\Q$, while
positivity of $a_i$ is equivalent to $\sigma=4h$.
The claim for signature $-4h$ can be deduced
by switching from $x$ to $-x$. \qed
\smallskip

Using Lemma \a.1 we can now give a complete description of
possible Witt classes of $SL(2,\Q)$-bundles.
\smallskip
\bf Theorem \a.\t. \rm(Theorem C)\sl\par
The set of Witt classes of all representations of $\pi_1(\Sigma_g)$
in $SL(2,\Q)$
is equal to the set of elements of $I^2\Q$ with norm
$\le4(g-1)$.
\smallskip\rm
Proof.
The classical Milnor--Wood inequality states that
the Euler class of a flat $SL(2,\R)$-bundle over 
$\Sigma_g$ has absolute value $\le g-1$. For an $SL(2,\Q)$-bundle
(treated as a flat $SL(2,\R)$-bundle)
this Euler class is equal to ${1\over 4}$ of the signature of 
the Witt class (cf.~[\TccI, Theorem 13.4]); therefore, the Witt class
of an $SL(2,\Q)$-bundle over $\Sigma_g$ has signature of absolute value
$\le4(g-1)$. Then it has also norm $\le4(g-1)$, by Lemma \a.1---except, possibly,
for $g=1$. But for $g=1$ we know that the Witt class is $0$ (by equicommutativity and Lemma 3.1),
so that the norm bound holds also in this case.

Conversely, by Theorem 11.6, every element of $I^2\Q$ 
of norm $\le 4(g-1)$
is realizable as the Witt class of some
$SL(2,\Q)$-bundle over $\Sigma_g$.\qed
\smallskip
\bf Remark \a.\t.\rm\par
Assume $g\ge2$. Then in the statement of Theorem \a.2 one can replace
``with norm $\le4(g-1)$'' by ``with signature of absolute value $\le4(g-1)$''
(as is evident from the proof).

\medskip
\bf 13. The easy norm bound is not sharp.\rm
\def\a{13}
\m=1
\n=1
\medskip
The Witt class $\w(\xi)$ of a flat $SL(2,K)$-bundle $\xi$ over $\Sigma_g$
has norm $\le 4g-2$.
For $K=\R$ and the Euler class this bound can be
improved to $4g-4$ (the already mentioned Milnor--Wood). It is unclear
whether this stronger estimate ($\|\w(\xi)\|\le 4g-4$)
holds for the general Witt class; we have not found any counterexamples.
The question is meaningful for fields with $u>4$. We now present an example
of an element of $I^2(K)$ with norm 6 that is not realizable as $\w(\xi)$
over $\Sigma_2$.
\smallskip
\bf Proposition \a.\t.\sl\par
Let $K=\Q((x))$, the field of Laurent series with rational coefficients.
Let
$$q=\langle 1,1,1,7,x,-7x\rangle.\leqno(\a.\f)$$
Then $q\in I^2(K)$, $\|q\|=6$ and $q$ is not
realizable as the Witt class of an $SL(2,K)$-bundle over $\Sigma_2$.
\rm
\smallskip
Proof.
We have $\dim(q)=6$ and $d_\pm(q)=1$, hence $q\in I^2(K)$.
Suppose that $q=\w(\xi)$. Let $\pi\colon\Sigma_3\to\Sigma_2$
be a degree-2 covering map. Then $\w(\pi^*\xi)=2\w(\xi)=2q$. Since
$\Sigma_3$ has a $\Delta$-complex structure with 10 triangles,
we would have $\|2q\|=\|\w(\pi^*\xi)\|\le10$. We will get a contradiction
by showing that
$\|2q\|=12$. (This equality also implies $\|q\|=6$.)

We will show that $q\oplus q$ is an anisotropic form.
To prove it we use (iteratively) the following fact (cf.~[Lam, Proposition VI.1.9]).
Let $F$ be a nondyadic complete discretely valued field with
residue field $\ov{F}$ and uniformizer $t$.
Then the form $q_1\oplus tq_2$ is anisotropic over $F$ if 
the forms $\ov{q}_1,\ov{q}_2$ over $\ov{F}$
are anisotropic.
We will use this
proposition for $F=\Q((x)),\Q_7$.

First,  we split
$$2q=\langle1,1,1,1,1,1,7,7\rangle\oplus x\langle 1,1,-7,-7\rangle.\leqno(\a.\f)$$
Thus, we want to show that
$\langle1,1,1,1,1,1,7,7\rangle$ and $\langle 1,1,-7,-7\rangle$
are anisotropic over $\Q$. 
For the first form this is clear---it is positive--definite.
We extend $\langle1,1,-7,-7\rangle$ to $\Q_7$ and split again:
$$\langle1,1,-7,-7\rangle=\langle1,1\rangle\oplus(-7)\langle1,1\rangle.\leqno(\a.\f)$$
Now $\langle 1,1\rangle$ is anisotropic over ${\bf F}_7$,
which finishes the proof.
\qed
\smallskip
Not all forms of norm 6 are susceptible to this argument, however.
For example, consider
$$q'=\langle1,1,1,5,x,-5x\rangle\leqno(\a.\f)$$
over the same field $\Q((x))$. This form is in $I^2\Q((x))$,
has norm 6, but $2q'$ is isotropic, hence has norm $\le10$.
We do not know whether $q'$ is realizable as the Witt class over $\Sigma_2$.

\bigskip
\bf References

\rm

\smallskip
\item{[Brown]}
 K.~Brown,
 Cohomology of groups.
 Graduate Texts in Mathematics, 87.
 Springer-Verlag, New York-Berlin, 1982.

\smallskip
\item{[\TccI]}
 J.~Dymara, T.~Januszkiewicz,
 Tautological characteristic classes I,
 2023.
 {\tt arXiv:2307.05765}.

\smallskip
\item{[EKM]}
 R.~Elman, N.~Karpenko and A.~Merkurjev,
 The algebraic and geometric theory of quadratic forms.
 American Mathematical Society Colloquium Publications, 56.
 American Mathematical Society, Providence, RI, 2008. 

\smallskip
\item{[GMS]}
 A.~Ghosh, C.~Meiri, P.~Sarnak, Commutators in ${\rm SL}_2$ and
 Markoff surfaces I. New Zealand J.~Math. {\bf 52} (2021 [2021--2022]),
 773--819. 

\smallskip
\item{[GSz]}
 P.~Gille and T.~Szamuely,
 Central simple algebras and Galois cohomology. 
 Cambridge Studies in Advanced Mathematics, 165.
 Cambridge University Press, Cambridge, 2017. 

\smallskip
\item{[Gold86]}
 W.~Goldman,
 Invariant functions on Lie groups and Hamiltonian flows of surface group representations.
 Invent.~Math. {\bf 85} (1986), no.~2, 263--302.

\smallskip
\item{[Gold88]}
 W.~Goldman,
 Topological components of spaces of representations.
 Invent.~Math. {\bf 93} (1988), no.~3, 557--607.

\smallskip
\item{[Hatcher]}
 A.~Hatcher, Algebraic Topology.
 Cambridge University Press, Cambridge, 2002.

\smallskip
\item{[Kr-T]}
 L.~Kramer, K.~Tent, A Maslov cocycle for unitary groups.
 Proc.~Lond.~Math.~Soc. (3) {\bf 100} (2010), no.~1, 91--115.

\smallskip
\item{[Lab]}
 F.~Labourie,
 Lectures on representations of surface groups.
 Zur.~Lect.~Adv.~Math.,
 European Mathematical Society (EMS), Z\"urich, 2013.

\smallskip
\item{[Lam]}
T.Y.~Lam,
Introduction to quadratic forms over fields.
Graduate Studies in Mathematics, 67. American Mathematical Society, Providence, RI, 2005.

\smallskip
\item{[Mats]}
 H.~Matsumoto,
 Sur les sous-groupes arithm\'etiques des groupes semi-simples d\'eploy\'es.
 Ann.~Sci.~\'Ecole Norm.~Sup.~(4) {\bf 2} (1969), 1--62.

\smallskip
\item{[Mil]}
 J.~Milnor,
 On the existence of a connection with curvature zero.
 Comment.~Math.~Helv. {\bf 32} 1958, 215--223.

\smallskip
\item{[MH]}
 J.~Milnor, D.~Husemoller,
 Symmetric bilinear forms. Ergebnisse der Mathematik und ihrer Grenzgebiete, Band 73. Springer--Verlag,
 New York--Heidelberg, 1973.

\smallskip
\item{[Moore]}
 C.~Moore, Group extensions of $p$-adic and adelic linear groups.
 Inst.~Hautes \'Etudes Sci.~Publ.~Math., {\bf 35} 1968, 157--222.

\smallskip
\item{[Ne]}
 J.~Nekov\'a$\check{\rm r}$,
 The Maslov index and Clifford algebras.
 (Russian) Funktsional.~Anal.~i Prilozhen. {\bf 24} (1990), no.~3, 36--44, 96;
 translation in Funct.~Anal.~Appl. {\bf 24} (1990), no.~3, 196--204 (1991). 

\smallskip
\item{[St]}
 R.~Steinberg, 
 Lectures on Chevalley groups. Notes prepared by John Faulkner and Robert Wilson.
 Revised and corrected edition of the 1968 original.
 University Lecture Series, 66. American Mathematical Society, Providence, RI, 2016. 

\smallskip
\item{[Su]}
 A.~Suslin,
 Torsion in $K_2$ of Fields,
 K-Theory {\bf 1} (1987), 5--29. 

\smallskip
\item{[Takeuchi]}
 K.~Takeuchi,
 Fuchsian groups contained in $SL_{2}(\Q)$,
 J.~Math.~Soc.~Japan, {\bf 23} (1971), no.~1, 82--94.

\smallskip
\item{[Wol]}
 S.~Wolpert,
 An elementary formula for the Fenchel--Nielsen twist.
 Comment.~Math. Helv. {\bf 56} (1981), no.~1, 132--135.
\bye